\documentclass[12pt]{article}
\usepackage{amssymb}
\usepackage{graphicx}
\usepackage{latexsym}
\def\part#1{\frac{\partial\phantom{q}}{\partial#1}}

\newenvironment{rmk}{\begin{trivlist}\item[]{\bf Remark:} }
{\end{trivlist}}
\newenvironment{rmks}{\begin{trivlist}\item[]{\bf Remarks:} }
{\end{trivlist}}
\newenvironment{ex}{\begin{trivlist}\item[]{\bf Example:} }
{\end{trivlist}}
\newenvironment{exs}{\begin{trivlist}\item[]{\bf Examples:} }
{\end{trivlist}}
\newenvironment{prf}{\begin{trivlist}\item[]{\bf Proof:} }
{\hfill $\Box$ \end{trivlist}}
\newenvironment{lemprf}{\begin{trivlist}\item[]{\bf Proof:} }
 {\end{trivlist}}
\newtheorem{thm}{Theorem}
\newtheorem{definition}{Definition}
\newtheorem{prp}[thm]{Proposition}
\newtheorem{lemma}[thm]{Lemma}

\newcommand{\lie}[1]{\mathfrak{#1}}
\def\End{\mathop{\rm End}\nolimits}
\def\Hom{\mathop{\rm Hom}\nolimits}

\def\deg{\mathop{\rm deg}\nolimits}
\def\adj{\mathop{\rm adj}\nolimits}

\def\Diff{\mathop{\rm Diff}\nolimits}
\def\Cliff{\mathop{\rm C l}\nolimits}

\newcommand{\R}{\mathbf{R}}
\newcommand{\C}{\mathbf{C}}

\newcommand{\Z}{\mathbf{Z}}

\newcommand{\CP}{{\mathbf C}{\rm P}}
\newcommand{\PP}{{\mathbf {\rm P}}}

\begin{document}
\title{Lectures on generalized geometry}
 \author{Nigel Hitchin\\[5pt]}
 \maketitle
\centerline{{\it Subject classification}: {Primary 53D18}}

\subsection*{Preface}

These notes are  based on six lectures given in March 2010 at the Institute of Mathematical Sciences in the  Chinese University of Hong Kong as part of the JCAS Lecture Series. They were mainly targeted at graduate students. They are not intended to be a comprehensive treatment of the subject of generalized geometry, but instead I have attempted to present the general features and to focus on a few topics which I have found particularly interesting and which I hope the reader will too. The relatively new material  consists of an account of Goto's existence theorem for generalized K\"ahler structures, examples of generalized holomorphic bundles and the B-field action on their moduli spaces.

 Since the publication of the first paper on the subject \cite{NJH1}, there have been many articles written within both  the mathematical and theoretical physics communities, and the reader should be warned that different authors have different conventions (or occasionally this author too!).  For other  accounts  of generalized geometry, I should direct the reader to the papers and surveys (e.g. \cite{Cav}, \cite{MG1}) of my former students Marco Gualtieri and Gil Cavalcanti who have developed many aspects of the theory. 

I would like to thank the IMS for its hospitality and for its invitation to give these  lectures, and Marco Gualtieri for useful conversations during the preparation of these  notes. 

 \tableofcontents

\section{The Courant bracket,  B-fields and metrics}

\subsection{Linear algebra preliminaries}
Generalized geometry is based on two premises -- the first is to replace the tangent bundle $T$ of a manifold $M$ by $T\oplus T^*$, and the second to replace the Lie bracket on sections of $T$ by the Courant bracket. The idea then is to use one's experience of differential geometry and by analogy to 
define and develop the generalized version. Depending on the object, this may or may not be a fertile process, but the intriguing fact is that, by drawing on the intuition of a mathematician, one  may often obtain this way a topic which is also of interest to the theoretical physicist. 

We begin with the natural linear algebra structure of the generalized tangent bundle  $T\oplus T^*$. 
If $X$ denotes a tangent vector and $\xi$  a cotangent vector then we write $X+\xi$ as a typical element of a fibre $(T\oplus T^*)_x$. There is a natural indefinite inner product defined by 
$$(X+\xi,X+\xi)=i_X\xi\quad (= \langle \xi,X\rangle =\xi(X))$$
using the interior product $i_X$, or equivalently the natural pairing $\langle \xi,X\rangle$ or the evaluation $\xi(X)$ of $\xi\in T^*_x$ on $X$. This is to be thought of as replacing the notion of a Riemannian metric, even though on an $n$-manifold it has signature $(n,n)$.

In block-diagonal form, a skew-adjoint transformation of $T\oplus T^*$ at a point can be written as 
$$\pmatrix {A & \beta\cr
                   B & -A^t}.$$
Here $A$ is just an endomorphism of $T$ and $B:T\rightarrow T^*$ to be skew-adjoint must satisfy
$$(B(X_1+\xi_1),X_2+\xi_2)= (B(X_1),X_2)=-(B(X_2),X_1)$$
so that $B$ is a skew-symmetric form, or equivalently $B\in \Lambda^2 T^*$, and its action is 
$X+\xi\mapsto i_XB.$
Since 
$$\pmatrix {0 & 0\cr
                   B & 0}^2= 0$$
  exponentiating gives 
  \begin{equation}
  X+\xi\mapsto X+\xi+i_XB
  \label{B}
  \end{equation}        
      This {\it B-field      action} will be fundamental, yielding extra transformations in generalized geometries.  It represents a breaking of symmetry in some sense since the bivector $\beta\in \Lambda^2T$ plays a lesser role. 
      
      \subsection{The Courant bracket}
      We described above the pointwise structure of the generalized tangent bundle. Now we consider the substitute for the Lie bracket $[X,Y]$ of two vector fields. This is the {\it Courant bracket} which appears in the literature in two different formats -- here we adopt the original skew-symmetric one: 
      
      \begin{definition} \label{cou} The Courant bracket of two sections $X+\xi,Y+\eta$ of $T\oplus T^*$ is defined by 
      $$[X+\xi,Y+\eta]=[X,Y]+{\mathcal L}_X\eta-{\mathcal L}_Y\xi-\frac{1}{2}d(i_X\eta-i_Y\xi).$$
      \end{definition}
      
      This has the important property that it commutes with the B-field action of a  {\it closed} 2-form $B$:
      
      \begin{prp} \label{Baction}
      Let $B$ be a closed 2-form, then
      $$[X+\xi+i_XB,Y+\eta+i_YB]=[X+\xi,Y+\eta]+i_{[X,Y]}B$$
      \end{prp}
      
      \begin{prf} We  shall make use of the Cartan formula for the Lie derivative of a differential form $\alpha$: ${\mathcal L}_X\alpha=d(i_X\alpha)+i_Xd\alpha$. First expand 
      $$[X+\xi+i_XB,Y+\eta+i_YB]=[X+\xi,Y+\eta]+{\mathcal L}_Xi_YB-{\mathcal L}_Yi_XB-\frac{1}{2}d(i_Xi_YB-i_Yi_XB).$$
      The last two terms give 
      $d(i_Yi_XB)={\mathcal L}_Yi_XB-i_Yd(i_XB)$ by the Cartan formula, and so yield 
     \begin{eqnarray*}
     [X+\xi+i_XB,Y+\eta+i_YB]&=&[X+\xi,Y+\eta]+{\mathcal L}_Xi_YB-i_Yd(i_XB)\\  
      &=& [X+\xi,Y+\eta]+i_{[X,Y]}B+i_Y{\mathcal L}_XB-i_Yd(i_XB)\\
       &=& [X+\xi,Y+\eta]+i_{[X,Y]}B+i_Yi_XdB
       \end{eqnarray*}
       by the Cartan formula again. So if $dB=0$ the bracket is preserved. 
      \end{prf}

    The inner product and Courant bracket naturally defined above   are clearly invariant under the induced action of a diffeomorphism of the manifold $M$. However, we now see that  a global closed differential 2-form $B$ will also act, preserving both the inner product and bracket. This means an overall action of the semi-direct product of closed 2-forms with diffeomorphisms
    $$\Omega^2(M)_{cl}\rtimes \Diff(M).$$
    This is a key feature of generalized geometry -- we have to consider B-field transformations as well as diffeomorphisms. 
    
   The Lie algebra of the group $\Omega^2(M)_{cl}\rtimes \Diff(M)$ consists of sections $X+B$ of $T\oplus \Lambda^2T^*$ where $B$ is closed.    If we take $B=-d\xi$, then the Lie algebra action on $Y+\eta$ is
   $$(X-d\xi)(Y+\eta)={\mathcal L}_X(Y+\eta)-i_Yd\xi=[X,Y]+{\mathcal L}_X\eta-{\mathcal L}_Y\xi+d(i_Y\xi).$$
 It is then easy to see  that we can reinterpret  the Courant bracket as the skew-symmetrization of this: 
   $$[X+\xi,Y+\eta]=\frac{1}{2}((X-d\xi)(Y+\eta)-(Y-d\eta)(X+\xi)).$$ 
   
However, although the Courant bracket is derived this way from a Lie algebra action, it is not itself a bracket of any Lie algebra -- the Jacobi identity fails. More precisely we have (writing $u=X+\xi, v=Y+\eta, w=Z+\zeta$)
\begin{prp}\label{jac}
$$[[u,v],w]+[[v,w],u]+[[w,u],v]=\frac{1}{3}d(([u,v],w)+([v,w],u)+([w,u],v))$$
\end{prp}
\begin{prf} If $u=X+\xi$, let $\tilde u=X-d\xi$ be the corresponding element in the Lie algebra of $\Omega^2(M)_{cl}\rtimes \Diff(M).$ We shall temporarily write $uv$ for the action of $\tilde u$ on $v$ (this is also called the Dorfman ``bracket" of $u$ and $v$) so that the Courant bracket is $(uv-vu)/2$. We first show that
\begin{equation}
u(vw)=(uv)w+v(uw).
\label{deriv}
\end{equation}
To see this note that $u(vw)-v(uw)=\tilde u\tilde v(w)-\tilde v\tilde u(w)=[\tilde u,\tilde v](w)$ since $\tilde u,\tilde v$ are Lie algebra actions, and the bracket here is just the commutator. But $(uv)w$ is the Lie algebra action of 
$uv=\tilde u v= [X,Y]+{\mathcal L}_X\eta-i_Yd\xi$
which acts as $[X,Y]-d({\mathcal L}_X\eta-i_Yd\xi)=[X,Y]-{\mathcal L}_Xd\eta+{\mathcal L}_Yd\xi$ using the Cartan formula and $d^2=0$. This however is just the  bracket $[\tilde u,\tilde v]$ in the Lie algebra of  
 $\Omega^2(M)_{cl}\rtimes \Diff(M).$
 
To prove the Proposition we note now  that the {\it symmetrization}  $(uv+vu)/2$ is
$$\frac{1}{2}({\mathcal L}_X\eta-i_Yd\xi+{\mathcal L}_Y\xi-i_Xd\eta)=\frac{1}{2}d(i_X\eta+i_Y\xi)=d(u,v)$$
while we have already seen that the skew-symmetrization $(uv-vu)/2$ is equal to $[u,v]$. So we rewrite the left hand side of the expression in the Proposition as one quarter of 
\begin{eqnarray*}
 (uv-vu)w - w(uv-vu)\\
            +(vw-wv)u- u(vw-wv)\\
            +(wu-uw)v-v(wu-uw)
          \end{eqnarray*}
          Using (\ref{deriv}) we sum these to get $(-1)$ times the sum $r$ of the right-hand column. If $\ell$ is the sum of the left-hand column then this means $\ell+r=-r$. But then $\ell-r=3(\ell+r)$ is the sum of terms like $(uv-vu)w+w(uv-vu)=4([u,v],w)$. The formula follows directly.
\end{prf}

There are two more characteristic properties of the Courant bracket which are easily verified:
\begin{equation}
[u,fv]=f[u,v]+(Xf)v-(u,v)df
\label{cour1}
\end{equation}
where $f$ is a smooth function, and as usual $u=X+\xi$, and 
\begin{equation}
X(v,w)=([u,v]+d(u,v),w)+(v,[u,w]+d(u,w)).
\label{cour2}
\end{equation}
\subsection{Riemannian geometry}
The fact that we introduced the inner product on $T\oplus T^*$ as the analogue of the Riemannian metric does not mean that Riemannian geometry is excluded from this area -- we just have to treat it in a different way.  We describe   a metric $g$ as a map $g:T\rightarrow T^*$ and consider its graph $V\subset T\oplus T^*$. This is the set of pairs $(X,gX)$ or in local coordinates (and the summation convention, which we shall use throughout)  the span of 
$$\frac{\partial}{\partial x_i}+g_{ij}dx_j.$$
The subbundle $V$ has an orthogonal complement $V^{\perp}$ consisting of elements of the form  $X-gX$. The inner product on $T\oplus T^*$ restricted to $X+gX\in V$ is $i_XgX=g(X,X)$ which is positive definite and restricted to $V^{\perp}$ we get the negative definite $-g(X,X)$. So $T\oplus T^*$ with its signature $(n,n)$ inner product can also be written as the orthogonal sum $V\oplus V^{\perp}$. Equivalently we have reduced the structure group of $T\oplus T^*$ from $SO(n,n)$ to $S(O(n)\times O(n))$.

The nondegeneracy of $g$ means that $g:T\rightarrow T^*$ is an isomorphism so that the projection from $V\subset T\oplus T^*$ to either factor is an isomorphism. This means we can lift vector fields or 1-forms to sections of $V$. Let us call $X^+$ the lift of a vector field $X$ to $V$ and $X^-$ its lift to $V^{\perp}$, i.e. $X^{\pm}=X\pm gX$. We also have the orthogonal projection $\pi_V:T\oplus T^*\rightarrow V$ and then 
$$\pi_V(X)=\pi_V\frac{1}{2}(X+gX + X-gX)=\frac{1}{2}X^+.$$

We can use these lifts and projections together with the Courant bracket to give a  convenient way of working out the Levi-Civita connection of $g$. 

First we show:
\begin{prp} \label{Vconnect} Let $v$ be a section of $V$ and $X$ a vector field, then 
$$\nabla_Xv=\pi_V[X^-,v]$$
defines a connection on $V$ which preserves the  inner product induced from $T\oplus T^*$.
\end{prp}
\begin{prf} Write $v=Y+\eta$, then observe that
$$\nabla_{fX}v=\pi_V[fX^-,v]=\pi_V(f[X^-,v]-(Yf)X^- +(v,X^-)df)$$
using Property (\ref{cour1}) of the Courant bracket. But $V$ and $V^{\perp}$ are orthogonal so $\pi_VX^-=0=(v,X^-)$ and hence $\nabla_{fX}v=f\nabla_Xv$.

Now using the same property we have
$$\nabla_{X}fv=\pi_V(f[X^-,v]+(Xf)v -(v,X^-)df)=f\nabla_{X}v+(Xf)v$$
 since $(v,X^-)=0$ and $\pi_Vv=v$. These two properties define a connection. 
 
 To show compatibility with the inner product take $v,w$ sections of $V$, then 
 $$(\nabla_Xv,w)+(v,\nabla_Xw)=(\pi_V[X^-,v],w)+(v,\pi_V[X^-,w])=([X^-,v],w)+(v,[X^-,w])$$
 since $\pi_V$ is the orthogonal projection onto $V$ and $v,w$ are sections of $V$. Now use  Property (\ref{cour2}) of the Courant bracket to see that 
 $$X(v,w)=([X^-,v]+d(X^-,v),w)+(v,[X^-,w]+d(X^-,w)).$$
 But $(X^-,v)=0=(X^-,w)$ and we get $X(v,w)=(\nabla_Xv,w)+(v,\nabla_Xw)$ as required.
\end{prf}

Using the isomorphism of $V$ with $T$ (or $T^*$) we can directly use this to find a connection on the tangent bundle. Directly, we take coordinates $x_i$ and then from the definition of the connection, the covariant derivative of ${\partial}/{\partial x_j}^+$ in the direction ${\partial}/{\partial x_i}$ is 
$$\pi_V\left[\frac{\partial}{\partial x_i}-g_{ik}dx_k, \frac{\partial}{\partial x_j}+g_{j\ell}dx_{\ell}\right].$$
Expanding the Courant bracket gives 
$$\frac{\partial g_{j\ell}}{\partial x_i}dx_{\ell}-\frac{\partial (-g_{ik})}{\partial x_j}dx_k-\frac{1}{2}d(g_{ji}-(-g_{ij}))=\frac{\partial g_{j\ell}}{\partial x_i}dx_{\ell}+\frac{\partial g_{ik}}{\partial x_j}dx_k-\frac{\partial g_{ij}}{\partial x_k}dx_{k}.$$
Projecting on $V$ we get
$$\frac{1}{2}(dx_k+g^{k\ell}\frac{\partial}{\partial x_{\ell}})\left(\frac{\partial g_{jk}}{\partial x_i}+\frac{\partial g_{ik}}{\partial x_j}-\frac{\partial g_{ij}}{\partial x_k}\right)=\frac{1}{2}g^{k\ell}\left(\frac{\partial g_{jk}}{\partial x_i}+\frac{\partial g_{ik}}{\partial x_j}-\frac{\partial g_{ij}}{\partial x_k}\right){\frac{\partial}{\partial x_{\ell}}}^{\!\!+}$$
which is the usual formula for the Christoffel symbols of the Levi-Civita connection.

\begin{ex}
Here is another computation -- the so-called Bianchi IX type metrics (using the terminology for example    in \cite{GP}). These are four-dimensional metrics with an $SU(2)$ action with generic orbit three-dimensional and in the diagonal form 
$$g=(abc)^2dt^2+a^2\sigma_1^2+b^2\sigma_2^2+c^2\sigma_3^2$$
where $a,b,c$ are functions of $t$ and $\sigma_i$ are basic left-invariant forms on the group, where $d\sigma_1=-\sigma_2\wedge \sigma_3$ etc. If $X_i$ are the dual vector fields then $[X_1,X_2]=X_3$  
and ${\mathcal L}_{X_1}\sigma_2=\sigma_3$ etc. 

Because of the even-handed treatment of forms and vector fields in generalized geometry, it is as easy to work out covariant derivatives of 1-forms as vector fields. Here we shall find the connection matrix for the  orthonormal basis of 1-forms $e_0=abc\, dt, e_1=a\sigma_1,e_2= b\sigma_2, e_3=c\sigma_3$. By symmetry it is enough to work out derivatives with respect to $X_1$ and $\partial/\partial t$. First we take $X_1$, so that $X_1^-=X_1-a^2\sigma_1$. 

For  the covariant derivative of $e_0$ consider the Courant bracket
$$[X_1-a^2\sigma_1,\frac{\partial}{\partial t} +  (abc)^2dt]=2aa'\sigma_1.$$
But $e^+_0=(abc)^{-1}(\partial/\partial t+(abc)^2dt)$ and using Property (\ref{cour1}) of the bracket and the orthogonality of $X_1^-,e_0^+$ we have 
$$[X_1^-,e_0^+]=\frac{2a'}{bc}\sigma_1=\frac{2a'}{abc}e_1$$
Projecting on $V$ and using $\pi_V e_1=e^+_1/2$, we have 
\begin{equation}
\nabla_{X_1}e_0=\frac{a'}{abc}e_1.
\label{10}
\end{equation}

For the 1-form $e_1$ note that, since ${\mathcal L}_{X_1}\sigma_1=0$
$$[X_1^-,X_1^+]=[X_1-a^2\sigma_1, X_1+a^2\sigma_1]=-\frac{1}{2}d(a^2+a^2)=-2aa'dt.$$
But $e_1=a^{-1}X_1^+$, and again using Property (\ref{cour1}) and the orthogonality of $X_1^-,X_1^+$ we have 
$[X_1^-,e_1^+]=-2a'dt$. Projecting onto $V$ gives $\pi_V(-2a' dt)= -(a'dt+a'(abc)^{-2}\partial/\partial t).$ So
\begin{equation}
\nabla_{X_1}e_1=-\frac{a'}{abc}e_0.
\label{11}
\end{equation}
(Note that with (\ref{10}) this checks with the fact that the connection preserves the metric.)
\vskip .25cm
For $e_2^+=b^{-1}X_2^+$ we have 
\begin{eqnarray*}
[X_1^-,X_2^+]&=&[X_1-a^2\sigma_1,X_2+b^2\sigma_2]=[X_1,X_2]+{\mathcal L}_{X_1}b^2\sigma_2+{\mathcal L}_{X_2}a^2\sigma_1-0\\
&=&X_3+(b^2-a^2)\sigma_3
\end{eqnarray*}
and so
$[X_1^-,e_2^+]=b^{-1}(X_3+(b^2-a^2)\sigma_3)$.
Projecting onto $V$,
$$\pi_V[X_1^-,e_2^+]=\frac{1}{2}b^{-1}(X_3+c^2\sigma_3)+\frac{1}{2}b^{-1}(b^2-a^2)(\sigma_3+c^{-2}X_3)$$
so that
\begin{equation}
\nabla_{X_1}e_2=\frac{1}{2bc}(c^2+b^2-a^2)e_3.
\label{12}
\end{equation}

Now we covariantly differentiate with respect to $t$.
\begin{eqnarray*}
\left[\frac{\partial}{\partial t}-(abc)^2dt, e_0^+\right]&=&\left[\frac{\partial}{\partial t}-(abc)^2dt, (abc)^{-1}\frac{\partial}{\partial t}+(abc) dt\right]\\
&=&-\frac{(abc)'}{(abc)^2}\frac{\partial}{\partial t}+(abc)'dt +(abc)'dt-\frac{1}{2}d(2(abc))\\
&=&\frac{(abc)'}{(abc)^2}\left(-\frac{\partial}{\partial t}+(abc)^2dt\right)
\end{eqnarray*}
so projecting onto $V$ gives 
\begin{equation}
\nabla_{\!\frac{\partial}{\partial t}}e_0=0.
\label{00}
\end{equation}
and finally  (similar to the first case above)
$$\left[\frac{\partial}{\partial t}-(abc)^2dt, X_1+a^2\sigma_1\right]=2aa'\sigma_1-0-0$$
so that
$$\left[\frac{\partial}{\partial t}-(abc)^2dt, e_1^+\right]=\left[\frac{\partial}{\partial t}-(abc)^2dt, a^{-1}(X_1+a^2\sigma_1)\right]=2a'\sigma_1-\frac{a'}{a^2}(X_1+a^2\sigma_1)$$
and projecting onto $V$ we get 
\begin{equation}
\nabla_{\!\frac{\partial}{\partial t}}e_1=\frac{a'}{a}e_1.
\label{01}
\end{equation}
\end{ex}

The point to make here is that the somewhat mysterious Courant bracket can be used as a tool for automatically computing covariant derivatives in ordinary Riemannian geometry.

 \section{Spinors, twists and skew torsion}
 \subsection{Spinors}\label{spin}
 In generalized geometry, the role of differential forms is changed. They become a {\it Clifford module} for the Clifford algebra generated by  $T\oplus T^*$ with its indefinite inner product. Recall that, given a vector space $W$ with an inner product $(\,\,,\,)$ the Clifford algebra $\Cliff(W)$ is generated by $1$ and $W$ with the relations $x^2=(x,x)1$ (in positive definite signature  the usual sign is $-1$ but this is the most  convenient for our case).
 
 Consider an exterior differential form $\varphi\in \Lambda^*T^*$ and define the action of $X+\xi\in T\oplus T^*$ on $\varphi$ by
 $$(X+\xi)\cdot \varphi=i_X\varphi+\xi\wedge \varphi$$
 then 
 $$(X+\xi)^2\cdot \varphi=i_X(\xi\wedge\varphi)+\xi\wedge i_X\varphi=i_X\xi\varphi=(X+\xi,X+\xi)\varphi$$
 and so  $\Lambda^*T^*$ is a module for the Clifford algebra. 
 
 We have already remarked that we can regard $T\oplus T^*$ as having structure group $SO(n,n)$ and if the manifold is oriented this lifts to $Spin(n,n)$. The exterior algebra is almost the basic spin representation of $Spin(n,n)$, but not quite. The Clifford algebra has an anti-involution -- any element is a sum of products $x_1x_2\dots x_k$ of generators $x_i\in W$ and 
 $$x_1x_2\dots x_k\mapsto x_kx_{k-1}\dots x_1$$
 defines the anti-involution. It represents a ``transpose" map $a\mapsto a^t$ arising from an invariant bilinear form on the basic spin module. In our case the spin representation is strictly speaking
 $$S=\Lambda^*T^*\otimes (\Lambda^nT^*)^{-1/2}.$$
 Another way of saying this is that there is an invariant  bilinear  form on $\Lambda^*T^*$ with values in the line bundle $\Lambda^nT^*$. Because of its appearance in another context it is known as the Mukai pairing. Concretely, given $\varphi_1,\varphi_2\in \Lambda^*T^*$, the pairing is 
 $$\langle \varphi_1,\varphi_2\rangle = \sum_j (-1)^j(\varphi_1^{2j}\wedge \varphi_2^{n-2j}+\varphi_1^{2j+1}\wedge \varphi_2^{n-2j-1})$$
 where the superscript $p$ denotes the $p$-form component of the form. 
 
 The Lie algebra of the spin group (which is the Lie algebra of $SO(n,n)$) sits inside the Clifford algebra as the subspace $\{a\in \Cliff(W): [a,W]\subseteq W\, {\mathrm {and}} \,\,  a=-a^t\}$ where the commutator is taken in the Clifford algebra. Consider  a 2-form $B\in \Lambda^2T^*$. The Clifford action of a 1-form $\xi$ is exterior multiplication, so $B\cdot \varphi=\sum b_{ij}\xi_i\cdot\xi_j\cdot\varphi=\sum b_{ij}\xi_i\wedge\xi_j\wedge\varphi$ defines an action on spinors.  Moreover, being skew-symmetric  in $\xi_i$ it satisfies $B^t=-B$. Now take  $X+\xi\in W=T\oplus T^*$ and the commutator $[B,X+\xi]$ in the Clifford algebra:
 $$B\wedge(i_X+\xi\wedge)\varphi-(i_X+\xi\wedge)B\wedge=B\wedge i_X\varphi-i_X(B\wedge\varphi)=-i_XB\wedge \varphi.$$
 So this action preserves $T\oplus T^*$ and so defines an element in the Lie algebra of $SO(n,n)$. But 
 the Lie algebra action of $B\in \Lambda^2T^*$ on $T\oplus T^*$ was $X+\xi\mapsto i_XB$, so we see from the above formula that the action of a B-field on  spinors is given by the exponentiation of $-B$ in the exterior algebra :
 $$\varphi\mapsto e^{-B\wedge}\varphi.$$
  One may easily check, for example,  that the Mukai pairing is invariant under the action: $\langle e^{-B}\varphi_1,e^{-B}\varphi_2\rangle=\langle\varphi_1,\varphi_2\rangle.$
  
  This action, together with the natural diffeomorphism action on forms, gives a combined action of the group $\Omega^2(M)_{cl}\rtimes \Diff(M)$ and a corresponding action of its Lie algebra. We earlier considered  the map $u\mapsto \tilde u$ given by $X+\xi\mapsto X-d\xi$ and on any bundle associated to $T\oplus T^*$  by a representation of $SO(n,n)$ we have an action of $\tilde u$. We regard this now as  a ``Lie derivative"  ${\mathbf L}_u$ in the direction of a section $u$ of $T\oplus T^*$. In the spin representation there is a ``Cartan formula" for this: 
  \begin{prp} The Lie derivative of a form $\varphi$ by a section $u$ of $T\oplus T^*$ is given by $${\mathbf L}_u\varphi=d(u\cdot \varphi)+u\cdot d\varphi.$$
  \end{prp}
 \begin{prf} $$d(X+\xi)\cdot \varphi+(X+\xi)\cdot d\varphi=di_X\varphi+d(\xi\wedge \varphi)+i_Xd\varphi+\xi\wedge d\varphi={\mathcal L}_X\varphi+d\xi\wedge\varphi$$
 using the usual Cartan formula and the fact that $B=-d\xi$ acts as $-B=d\xi$.
 \end{prf}
In fact  replacing the exterior product by the Clifford product is a common  feature of generalized geometry whenever we deal with forms.

The Lie derivative acting on sections of $T\oplus T^*$ is the Lie algebra action we observed in the first lecture so 
\begin{equation} 
{\mathbf L}_uv-{\mathbf L}_vu=2[u,v]
\label{Lie}
\end{equation}
where $[u,v]$ is the Courant bracket.

\subsection{Twisted structures}\label{twist}
We now want to consider a twisted version of $T\oplus T^*$. Suppose we have a nice covering of the manifold $M$ by open sets $U_{\alpha}$ and we give ourselves a closed 2-form $B_{\alpha\beta}=-B_{\beta\alpha}$ on each two-fold intersection $U_{\alpha}\cap U_{\beta}$. We can use the action of  $B_{\alpha\beta}$  to identify $T\oplus T^*$ on $U_{\alpha}$ with $T\oplus T^*$ on $U_{\beta}$ over the intersection. This will be compatible over threefold intersections if 
\begin{equation}
B_{\alpha\beta}+B_{\beta\gamma}+B_{\gamma\alpha}=0
\label{cocycle}
\end{equation}
on $U_{\alpha}\cap U_{\beta}\cap U_{\gamma}$. 

We have seen in the first lecture that the action of a closed 2-form on $T\oplus T^*$ preserves both the inner product and the Courant bracket, so by the above identifications this way we construct a rank $2n$ vector bundle $E$ over $M$ with an inner product and a  bracket operation on sections. And since the B-field action is trivial on $T^*\subset T\oplus T^*$, the vector bundle is an extension:
$$0\rightarrow T^*\rightarrow E\stackrel{\pi}\rightarrow T\rightarrow 0.$$
Such an object is called an {\it exact Courant algebroid}. It can be abstractly characterized by the Properties (\ref{cour1}) and (\ref{cour2}) of the Courant bracket, where the vector field $X$ is $\pi u$, together with the Jacobi-type formula in Proposition \ref{jac}. 

 The relation (\ref{cocycle}) says that we have a 1-cocycle for the sheaf $\underline\Omega^2_{cl}$ of closed 2-forms on $M$. There is an exact sequence of sheaves 
 $$0\rightarrow \underline\Omega^2_{cl}\rightarrow \underline\Omega^2\stackrel{d}\rightarrow \underline\Omega^3_{cl}\rightarrow 0$$
and since $\underline\Omega^2$ is a flabby sheaf, we have 
$$H^1(M,\underline\Omega^2_{cl})\cong H^0(M,{\underline\Omega^3_{cl}})/dH^0(M,{\underline\Omega^2})=\Omega_{cl}^3/d\Omega^2= H^3(M,\R)$$
so that such a structure has a characteristic degree $3$ cohomology class.

\begin{ex} The theory of gerbes fits into the twisted picture quite readily. Very briefly, a $U(1)$ gerbe can be defined by a 2-cocycle with values on the sheaf of $C^{\infty}$ circle-valued functions -- so it is given by functions $g_{\alpha\beta\gamma}$ on threefold intersections satisfying a coboundary condition. In the exact sequence of sheaves of $C^{\infty}$ functions
$$1\mapsto \Z\mapsto \underline\R\stackrel{exp \,2\pi i}\rightarrow \underline U(1)\rightarrow 1$$
the 2-cocycle defines a class in $H^3(M,\Z)$.

If we think of the analogue for line bundles, we have the transition functions $g_{\alpha\beta}$ and then a connection on the line bundle is given by 1-forms $A_{\alpha}$ on open sets such that
$$A_{\beta}-A_{\alpha}=(g^{-1}dg)_{\alpha\beta}$$
on twofold intersections (where we identify the Lie algebra of the circle with $\R$).

A {\it connective structure} on a gerbe is similarly a collection of 1-forms $A_{\alpha\beta}$ such that 
$$A_{\alpha\beta}+A_{\beta\gamma}+A_{\gamma\alpha}=(g^{-1}dg)_{\alpha\beta\gamma}$$
on threefold intersections. Clearly $B_{\alpha\beta}=dA_{\alpha\beta}$ defines a Courant algebroid, and its characteristic class is the image of the integral cohomology class in $H^3(M,\R)$.

An example of this is a hermitian structure on a holomorphic gerbe on a complex manifold (defined by a cocycle of holomorphic functions $h_{\alpha\beta\gamma}$ with values in $\C^*$). A hermitian structure on this is a choice of a cochain $k_{\alpha\beta}$ of positive functions  with $\vert h_{\alpha\beta\gamma}\vert =k_{\alpha\beta}k_{\beta\gamma}k_{\gamma\alpha}$. Then $A_{\alpha\beta}=(h^{-1}d^ch)_{\alpha\beta}$ defines a connective structure. Here $d^c=I^{-1}dI=-i(\partial-\bar\partial)$.
\end{ex}
\vskip .25cm
The bundle $E$ has an orthogonal structure and so an associated spinor bundle $S$. By the definition of $E$, $S$ is obtained by identifying $\Lambda^*T^*$ over $U_{\alpha}$ with $\Lambda^*T^*$ over $U_{\beta}$ by 
$$\varphi\mapsto e^{-B_{\alpha\beta}}\varphi.$$
A global section of  $S$ is then given by local forms $\varphi_{\alpha},\varphi_{\beta}$ such that  $\varphi_{\alpha}= e^{-B_{\alpha\beta}}\varphi_{\beta}$ on $U_{\alpha}\cap U_{\beta}$. Since $B_{\alpha\beta}$ is closed and even,
$$d \varphi_{\alpha}= d(e^{-B_{\alpha\beta}}\varphi_{\beta})=-dB_{\alpha\beta}\wedge e^{-B_{\alpha\beta}}\varphi_{\beta}+e^{-B_{\alpha\beta}}d\varphi_{\beta}=e^{-B_{\alpha\beta}}d\varphi_{\beta}$$
and so we have a well-defined operator
$$d:C^{\infty}(S^{ev})\rightarrow C^{\infty}(S^{od}).$$
The $\Z_2$-graded cohomology of this is the {\it twisted cohomology}. There is a more familiar way of writing this if we consider the inclusion of sheaves $\underline\Omega^2_{cl}\subset \underline\Omega^2$. Since $\underline\Omega^2$ is a flabby sheaf, the cohomology class of $B_{\alpha\beta}$ is trivial here and we can find 2-forms $F_{\alpha}$ such that on  $U_{\alpha}\cap U_{\beta}$
$$F_{\beta}-F_{\alpha}=B_{\alpha\beta}.$$
Since $B_{\alpha\beta}$ is closed  $dF_{\beta}=dF_{\alpha}$ is the restriction of a global closed 3-form $H$ which represents the characteristic class in $H^3(M,\R)$.

But then  
$$e^{-F_{\alpha}}\varphi_{\alpha}=e^{-F_{\beta}}e^{B_{\alpha\beta}}\varphi_{\alpha}=e^{-F_{\beta}}\varphi_{\beta}$$
defines a global exterior form $\psi$. Furthermore
$$d\psi=d(e^{-F_{\alpha}}\varphi_{\alpha})=-H\wedge \psi+e^{-F_{\alpha}}d\varphi_{\alpha}.$$
Thus the operator $d$ above defined on $S$ is equivalent to the operator
$$d+H:\Omega^{ev}\rightarrow \Omega^{od}$$
on exterior forms.

\begin{ex} For gerbes the full analogue of a connection is a connective structure together with a {\it curving}, which is precisely a choice of 2-form $F_{\alpha}$ such that $F_{\beta}-F_{\alpha}=dA_{\alpha\beta}$. In this case the 3-form $H$ such that $H/2\pi$ has integral periods is the curvature.
\end{ex} 

\begin{rmk} Instead of thinking in cohomological terms about writing a cocycle of closed 2-forms $B_{\alpha\beta}$  as a coboundary $F_{\beta}-F_{\alpha}$ in the sheaf of all 2-forms, there is a more geometric interpretation of this choice which  can be quite convenient. The B-field action of $F_{\alpha}$ gives an isomorphism of $T\oplus T^*$ with itself over $U_{\alpha}$ and the relation $F_{\beta}-F_{\alpha}=B_{\alpha\beta}$ says that this extends to an isomorphism
$$E\cong T\oplus T^*.$$
More concretely, $X$ over $U_{\alpha}$ is mapped to  $X+i_XF_{\alpha}\in E$ and defines a splitting (in fact an {\it isotropic} splitting)  of the extension $0\rightarrow T^*\rightarrow E\rightarrow T\rightarrow 0$. 

With the same coboundary data, we identified the spinor bundle $S$ with the exterior algebra bundle and now we note that
$$(X+i_XF_{\alpha})\cdot e^{-F_{\alpha}}\varphi_{\alpha}=-e^{-F_{\alpha}}i_X\varphi_{\alpha}$$
so the two are compatible. We have a choice -- either consider $E,S$ with their standard local models of $T\oplus T^*$ and $\Lambda^*T^*$, or make the splitting and give a global isomorphism. The cost is that we replace the ordinary exterior derivative by $d+H$ and, as can be seen from the proof of Proposition \ref{Baction}, replace the standard Courant bracket by the twisted version 
\begin{equation}
[X+\xi,Y+\eta]+i_Yi_XH.
\label{Couranttwist}
\end{equation}
\end{rmk}

\subsection{Skew torsion}\label{skewsection}
If we replace $T\oplus T^*$ by its twisted version $E$ we may ask how to incorporate a Riemannian metric as we did in the first lecture. Here is the definition:
\begin{definition} \label{geng} A generalized metric is a subbundle $V\subset E$ of rank $n$ on which the induced inner product is positive definite. 
\end{definition}

Since the inner product on $T^*\subset E$ is zero and is positive definite on $V$, $V\cap T^*=0$ and so in a local isomorphism $E\cong T\oplus T^*$, $V$ is the graph of a map $h_{\alpha}:T\rightarrow T^*$. So, splitting into symmetric and skew symmetric parts
$$h_{\alpha}=g_{\alpha}+F_{\alpha}.$$
On the twofold intersection
$$h_{\alpha}(X)=h_{\beta}(X)+i_XB_{\alpha\beta}.$$
Thus $h_{\alpha}(X)(X)=h_{\beta}(X)(X)=g(X,X)$ for a well-defined Riemannian metric $g$, but 
$F_{\alpha}=F_{\beta}+B_{\alpha\beta}$. Associated with a generalized metric we thus obtain a natural splitting.

Now the definition of a connection in Proposition \ref{Vconnect} makes perfectly good sense in the twisted case. To see what we get, let us redo the calculation using local coordinates. In this case $V$ is defined locally  by $X+gX+i_XF_{\alpha}$ and $V^{\perp}$ by  $X-gX+i_XF_{\alpha}$ (we are just transforming the $V$ and $V^{\perp}$ of the metric $g$ by the orthogonal transformation of the local B-field $F_{\alpha}$.) So the appropriate Courant bracket is 
$$\left[\frac{\partial}{\partial x_i}-g_{ik}dx_k+F_{ik}dx_k, \frac{\partial}{\partial x_j}+g_{j\ell}dx_{\ell}+F_{j\ell}dx_{\ell}\right]$$
and this gives the terms for the Levi-Civita connection plus a term
\footnote[1]{In a parallel discussion in \cite{NJH2} the third term in this expansion was unfortunately omitted.}

$$\frac{\partial F_{j\ell}}{\partial x_i}dx_{\ell}-\frac{\partial F_{ik}}{\partial x_j}dx_k-\frac{1}{2}d(F_{ji}-F_{ij})=\left(\frac{\partial F_{j\ell}}{\partial x_i}-\frac{\partial F_{i\ell}}{\partial x_j}+\frac{\partial F_{ij}}{\partial x_{\ell}}\right)dx_{\ell}.$$
(We could also have used the twisted bracket as in (\ref{Couranttwist}).)

Writing the skew bilinear form $F_{\alpha}$ as $\sum_{i<j}F_{ij}dx_i\wedge dx_j$, and $dF_{\alpha}=\sum_{i<j<k}H_{ijk}dx_i\wedge dx_j\wedge dx_k$, this last term is $H_{ji\ell}dx_{\ell}$ and represents a connection with {\it skew torsion}:  recall that the torsion of a connection on the tangent bundle is  
$$T(X,Y)=\nabla_XY-\nabla_YX-[X,Y]$$
and is said be skew  if $g(T(X,Y),Z)$ is skew-symmetric.
In our case, with $X=\partial/\partial x_i,Y=\partial/\partial x_j$, we use as before the projection onto $V$ to get
$$T\left(\frac{\partial}{\partial x_i},\frac{\partial}{\partial x_j}\right)=H_{ji\ell}g^{\ell k}\frac{\partial }{\partial x_k}.$$

\begin{rmk} We could have interchanged the roles of  $V$ and $V^{\perp}$ in the above argument -- this would give a connection with opposite torsion $-H$.
\end{rmk}

\begin{exs}

\noindent 1. The standard example of a connection with skew torsion is the flat connection given by trivializing the tangent bundle on a compact Lie group by left translation, i.e. for all $X\in {\lie{g}}$, $\nabla X=0$. Then 
$T(X,Y)=[X,Y]$.
 
\noindent 2. The second class of examples is given by the Bismut connection on the tangent bundle of a Hermitian manifold. In \cite{B}  it is shown that any Hermitian manifold has a unique connection with skew torsion which preserves the complex structure and the Hermitian metric. In the K\"ahler case it is the Levi-Civita connection, but in general the skew torsion is defined by the 3-form $d^c\omega$ where $\omega$ is the Hermitian form.
\end{exs}

\section{Generalized complex manifolds}
\subsection{Generalized complex structures}\label{gcs}
There are many ways of defining a complex manifold other than by the existence of holomorphic coordinates. One is the following: an endomorphism $J$ of $T$ such that $J^2=-1$ and such that the $+i$ eigenspace of $J$ on the complexification $T\otimes \C$ is Frobenius-integrable. Of course it needs the Newlander-Nirenberg theorem to find the holomorphic coordinates from this definition, but it is one that we can adapt straightforwardly to the generalized case. The only extra condition is compatibility with the inner product. So we have:
\begin{definition} A generalized complex structure on a manifold is an endomorphism $J$ of $T\oplus T^*$ such that
\begin{itemize}
\item
$J^2=-1$
\item
$(Ju,v)=-(u,Jv)$
\item
Sections of the subbundle $E^{1,0}\subset (T\oplus T^*)\otimes \C$ defined by the $+i$ eigenspaces of $J$ are closed under the Courant bracket.
\end{itemize}
\end{definition}

\begin{rmks} 

\noindent 1. It is not immediately obvious that the obstruction to integrability is tensorial, i.e. if $[u,v]$ is a section of $E^{1,0}$ then so is $[u,fv]$, but this is indeed so. It depends on the fact that $E^{1,0}$ is isotropic with respect to the inner product. In fact if $Ju=iu$ then 
$$i(u,u)=(Ju,u)=-(u,Ju)=-i(u,u).$$
Then recall Property (\ref{cour1}) of the Courant bracket: $[u,fv]=f[u,v]+(Xf)v-(u,v)df$. If $u,v$ are sections of $E^{1,0}$ then $(u,v)=0$ so 
$$[u,fv]=f[u,v]+(Xf)v$$
and if $u,v$ and $[u,v]$ are sections, so is $[u,fv]$.

\noindent 2. The definition obviously extends to the twisted case, replacing $T\oplus T^*$ by $E$.

\noindent 3. The data of $J$ is equivalent to giving an isotropic subbundle $E^{1,0}\subset (T\oplus T^*)\otimes \C$ of rank $n$ such that $E^{1,0}\cap \bar E^{1,0}=0$. This is the way we shall describe examples below, and we shall write $E^{0,1}$ for $\bar E^{1,0}$.

\noindent 4. The endomorphism $J$ reduces the structure group of $T\oplus T^*$ from $SO(2m,2m)$ to the indefinite unitary group $U(m,m)$.
\end{rmks}

\begin{exs}

\noindent 1. An ordinary complex manifold is an example. We take $E^{1,0}$ to be spanned by $(0,1)$ tangent vectors and $(1,0)$ forms:
$$\frac{\partial}{\partial \bar z_1}, \frac{\partial}{\partial \bar z_2},\dots , dz_1, dz_2,\dots$$
The integrability condition is obvious here.

\noindent 2. A symplectic form $\omega$ defines a generalized complex structure. Here $E^{1,0}$ is spanned by sections of $(T\oplus T^*)\otimes \C$ of the form
$$\frac{\partial}{\partial x_j}-i\omega_{jk}dx_k.$$
This is best seen as the transform of $T\subset T\oplus T^*$ by the complex B-field $-i\omega$. Since $T$ is isotropic and its sections are obviously closed under the Courant bracket (which is just the Lie bracket on vector fields)  the same is true of its transform by a closed 2-form.

\noindent 3. A holomorphic Poisson manifold is an example. Recall that a Poisson structure on a manifold is a section $\sigma$ of $\Lambda^2T$, which therefore defines a homomorphism $\sigma:T^*\rightarrow T$. For a function $f$, $\sigma(df)=X$ is   called a Hamiltonian vector field and $X(g)=\sigma(df,dg)$ is defined to be the Poisson bracket $\{f,g\}$ of the two functions. The integrability condition for a Poisson structure is
$$\sigma(d\{f,g\})=[X,Y]$$
where $Y$ is the Hamiltonian vector field of $g$. 

If $M$ is a complex manifold and $\sigma\in \Lambda^2T^{1,0}$ a holomorphic Poisson structure then $E^{1,0}$ is spanned by 
$$\frac{\partial}{\partial \bar z_1}, \frac{\partial}{\partial  \bar z_2},\dots,  dz_1-\sigma(dz_1), dz_2-\sigma(dz_1),\dots$$
Because $\sigma$ is holomorphic the only potentially non-trivial Courant brackets are of the form
$$[dz_i-\sigma(dz_i),dz_j-\sigma(dz_j)]=[\sigma(dz_i),\sigma(dz_j)]-d\{z_i,z_j\}+d\{z_j,z_i\}+\frac{1}{2}d(\{z_i,z_j\}-\{z_j,z_i\}).$$
But the Courant bracket on vector fields is the Lie bracket so by integrability of the Poisson structure 
$[\sigma(dz_i),\sigma(dz_j)]=\sigma(d\{z_i,z_j\})$ and hence the bracket above is 
$$\sigma(d\{z_i,z_j\})-d\{z_i,z_j\}$$
which again lies in $E^{1,0}$.
\end{exs}

There is another way to describe the integrability which can be very useful. The subbundle $E^{1,0}$ has rank $n$ and is isotropic in a $2n$-dimensional space. For a non-degenerate inner product this is the maximal dimension.  Given any spinor $\psi$, the space of $x\in W$ such that $x\cdot\psi=0$ is isotropic because $0=x\cdot x\cdot\psi=(x,x)\psi$. A maximal isotropic subspace is determined by a special type of spinor called a pure spinor. 

So to any maximal isotropic subspace we can associate a one-dimensional space of pure spinors it annihilates. Hence a generalized complex manifold has a complex line subbundle of $\Lambda^*T^*\otimes \C$ (called the {\it canonical bundle}) consisting of multiples of a pure spinor defining $J$. The condition $ E^{1,0}\cap \bar E^{1,0}=0$ is equivalent to the Mukai pairing $\langle \psi,\bar\psi\rangle$ for the spinor and its conjugate being non-zero. 

\begin{exs}

\noindent 1. For an  ordinary complex manifold the subspace $$\frac{\partial}{\partial \bar z_1}, \frac{\partial}{\partial \bar z_2},\dots , dz_1, dz_2,\dots$$
annihilates $dz_1\wedge dz_2\wedge\dots\wedge dz_m$. This generates the usual canonical bundle of a complex manifold.

\noindent 2. The tangent space $T$ annihilates $1$ by Clifford multiplication (interior product). Hence for a symplectic manifold the transform of $T$ by $-i\omega$ annihilates the form $e^{i\omega}$. Here the canonical bundle is trivialized by this form. 
\end{exs}

Here is integrability in this context:
\begin{prp} \label{int} Let $\psi$ be a form which is a pure spinor with $\langle \psi,\bar\psi\rangle\ne 0$. Then it defines a generalized complex structure if and only if $d\psi=w\cdot\psi$ for some local section $w$ of $(T\oplus T^*)\otimes \C$.
\end{prp}

\begin{prf} First assume $d\psi=w\cdot\psi$. Suppose $u\cdot \psi=0=v\cdot \psi$. Then, since the Lie derivative ${\mathbf L}_v$ acts via the Lie algebra action and so preserves the Clifford product, we have  
$0={\mathbf L}_v(u\cdot\psi)={\mathbf L}_vu\cdot\psi+u\cdot {\mathbf L}_v\psi.$

Using the Cartan formula ${\mathbf L}_v\psi=d(v\cdot\psi)+v\cdot d\psi=v\cdot d\psi=v\cdot w\cdot \psi$ 
and so
$${\mathbf L}_vu\cdot\psi+u\cdot v\cdot w\cdot\psi=0$$
Now use  the Clifford relations, 
$$u\cdot v\cdot w\cdot\psi=u\cdot(2(v,w)-w\cdot v)\cdot\psi=0$$
since  $u\cdot \psi=0=v\cdot \psi$. We deduce that ${\mathbf L}_vu\cdot\psi=0$.
Hence, interchanging the roles of $u$ and $v$ and subtracting, 
$$0={\mathbf L}_vu\cdot\psi-{\mathbf L}_uv\cdot\psi=2[v,u]\cdot \psi$$
from (\ref{Lie}). 

The Courant bracket therefore preserves the annihilator of $\psi$ and we have the integrability condition for a generalized complex structure.
\vskip .25cm
Conversely, assume the structure is integrable. If $u\cdot \psi=0=v\cdot \psi$ then from the definition of integrability  $[u,v]\cdot\psi=0$ and so
${\mathbf L}_vu\cdot\psi-{\mathbf L}_uv\cdot\psi=0$. But then from the above algebra we have
$(u\cdot v-v\cdot u)\cdot d\psi=0$ or, since $u\cdot v=- v\cdot u+2(u,v)1=- v\cdot u$,
$$u\cdot v\cdot d\psi=0.$$
As far as the linear algebra is concerned, any two endomorphisms $J$ satisfying the conditions of a generalized complex structure are equivalent under the action of $SO(2m,2m)$ -- they form the orbit $SO(2m,2m)/U(m,m)$. Hence to proceed, we can use the linear algebra of the standard complex structure where $\psi=dz_1\wedge\dots\wedge dz_m$ to determine those  $\varphi$ which satisfy $u\cdot v\cdot \varphi=0.$

 Taking $u=dz_i$ and $v=\partial/\partial \bar z_j$, the condition that $u\cdot v\cdot \varphi=0$ for all $i$ means that $\varphi$ is a sum of forms of type $(m,q)$ or $(p,0)$. Taking $u=dz_i,v=dz_j$ we must have $p=m$ or $m-1$. Taking $u=\partial/\partial \bar z_i, v=\partial/\partial \bar z_j$ we need $q=0$ or $1$. Thus $\varphi$ is a sum of $(m,0),(m,1)$ and $(m-1,0)$ terms. But $d\psi$ has opposite parity to $\psi$ so it must be $(m,1)$ and $(m-1,0)$ terms. However, these are generated by $d\bar z_i\wedge dz_1\wedge\dots\wedge dz_m$ and $i_{\partial/\partial  z_j}dz_1\wedge\dots\wedge dz_m$, that is $w\cdot \psi$, as required.
\end{prf}

\begin{ex} The simplest use of this integrability is when there is a global closed form which is a pure spinor. Such manifolds are called {\it generalized Calabi-Yau manifolds} and include ordinary Calabi-Yau manifolds where the holomorphic $m$-form is $\psi$, or symplectic manifolds where $\psi=e^{i\omega}$.

\end{ex}

\subsection{Symmetries and twisting}
At first sight, there seems little common ground when we think of the symmetries of symplectic manifolds or complex manifolds. In the first case, any smooth function defines a Hamiltonian vector field, in the second the Lie algebra of holomorphic vector fields is at most finite-dimensional and often zero.  In generalized geometry, however, we use the extended group $\Omega^2(M)_{cl}\rtimes \Diff(M)$ and this restores the balance between the two. 

\begin{prp} Let $J$ be a generalized complex structure and $f$ a smooth function. Then if $X+\xi=J(df)$, $X-d\xi$ in the Lie algebra of $\Omega^2(M)_{cl}\rtimes \Diff(M)$ preserves $J$.
\end{prp}
\begin{prf}
If $u=Jdf$, then decompose $u=u^{1,0}+u^{0,1}$ into its $\pm i$ eigenspace components of $J$, so that $-df=Ju= iu^{1,0}-iu^{0,1}$. Let $\psi$ be a local section of the canonical bundle, then $u^{1,0}\cdot\psi=0$ and hence
$$u\cdot\psi=u^{0,1}\cdot\psi=-idf\cdot\psi=-idf\wedge\psi.$$
Thus, using Proposition \ref{int} and the Cartan formula 
$${\mathbf L}_u\psi=d(u\cdot \psi)+u\cdot d\psi=d(-idf\wedge\psi)+u\cdot d\psi=(idf+u)\cdot w\cdot \psi.$$
But $idf+u=u^{1,0}-u^{0,1}+u^{1,0}+u^{0,1}=2u^{1,0}$ and so
$${\mathbf L}_u\psi=2u^{1,0}\cdot w\cdot\psi=4(u^{1,0},w)\psi$$
using the Clifford identity $u^{1,0}\cdot w+w\cdot u^{1,0}=2(u^{1,0},w)1$ and $u^{1,0}\cdot\psi=0$.

It follows that the Lie derivative of $\psi$ is a multiple of $\psi$ and so preserves  the generalized complex structure. 
\end{prf}

From this proposition we can see how a complex manifold acquires symmetries from smooth functions, for in this case $Jdf=X+\xi=-d^cf$ and exponentiating  we have the B-field action of the closed 2-form $dd^cf$.
\vskip .25cm
Slightly more generally, any real closed $(1,1)$-form is a symmetry of an ordinary complex structure thought of as a generalized complex structure, since the interior product with a $(0,1)$-vector $\partial/\partial \bar z_i$ is a $(1,0)$ form -- a linear combination of $dz_j$s. It follows that, if we take a 1-cocycle of such forms and construct an extension $E$ as in Section \ref{twist}, we obtain a twisted generalized complex structure on $E$. 

In particular, given a closed 3-form $H$ of type $(1,2)$ we can write this on a small enough  open set $U_{\alpha}$ as $\partial\bar\partial A_{\alpha}$ for a $(0,1)$-form $A_{\alpha}$. Then on a twofold intersection
$\partial\bar\partial A_{\alpha}-\partial\bar\partial A_{\beta}=0$ and hence $d(\bar\partial A_{\alpha}-\bar\partial A_{\beta})=0$. Defining $B_{\alpha\beta}$ to be the real part of $\bar\partial (A_{\beta}- A_{\alpha})$ gives such a cocycle.

\subsection{The $\bar\partial$-operator}\label{dbarsec}
On a  manifold with generalized complex structure $J$ we have the eigenspace decomposition $(T\oplus T^*)\otimes \C= E^{1,0}\oplus E^{0,1}$, and given a function $f$ we define $\bar\partial_Jf$ to be the $(0,1)$-component. Note that in the twisted case we have $T^*\subset E$ and so we can do the same. 
\begin{exs}

\noindent 1. For an ordinary complex manifold 
$$\bar\partial_Jf=\bar\partial f.$$

\noindent 2. For a symplectic structure $\omega$ 
$$\bar\partial_Jf=\frac{1}{2}(iX+df)$$
where $X$ is the Hamiltonian vector field of $f$ i.e.  $i_X\omega=df$. Note that in this case $\bar\partial_Jf=0$ implies that $f$ is constant, so we can't approach generalized complex geometry purely in terms of sheaves of local holomorphic functions. 

\noindent 3. For a holomorphic Poisson structure $\sigma$
$$\bar\partial_Jf=\bar\partial f+\sigma(\partial f)-\bar\sigma(\bar\partial f)$$
\end{exs}

The $\bar\partial$-operator maps functions to sections of $E^{0,1}$, and we want to extend it to a complex just like the Dolbeault complex of a complex manifold. For this we need to extend to an operator
$$\bar\partial_J:C^{\infty}(E^{0,1})\rightarrow C^{\infty}(\Lambda^2 E^{0,1}).$$
There is an obvious formula, analogous to the usual definition of  the exterior derivative, but using the Courant bracket instead of the Lie bracket.  Note first that because $E^{1,0}$ is maximal isotropic, the inner product identifies $E^{0,1}$ with the dual of $E^{1,0}$, so $\Lambda^p E^{0,1}$ is the space of alternating multilinear $p$-forms on $E^{1,0}$.

The extended operator is then given, for sections $u=X+\xi,v=Y+\eta$ of  $E^{1,0}$ by the usual formula, but with the Courant bracket $[u,v]$: 
\begin{equation}
2\bar\partial_J \alpha(u,v)=X(\alpha(v))-Y(\alpha(u))+\alpha([u,v]).
\label{da}
\end{equation}
Given this, there is an obvious extension to an operator
$$\bar\partial_J:C^{\infty}(\Lambda^pE^{0,1})\rightarrow C^{\infty}(\Lambda^{p+1} E^{0,1})$$
which satisfies the property
$$\bar\partial_J(f\alpha)=\bar\partial_Jf\wedge \alpha+f\bar\partial_J\alpha.$$

The only point to make here is that the relation $\bar\partial_J^2=0$ requires the Jacobi identity for the Courant bracket which we know doesn't hold in general. However, its failure as in Proposition \ref{jac} is due to the term
$$d(([u,v],w)+([v,w],u)+([w,u],v)).$$
But by the definition of a generalized complex structure the Courant bracket $[u,v]$ of two sections of $E^{1,0}$ is also a section, and as we have seen,  $E^{1,0}$ is isotropic. It follows that this expression is zero for such sections and the Jacobi identity holds. 

\begin{exs} 

\noindent 1. In the case of an ordinary complex structure $E^{0,1}=\bar T^*\oplus T$ where $T$ here denotes the holomorphic tangent bundle. Hence
$$\Lambda^m E^{0,1}=\bigoplus_{p+q=m}\Lambda^pT\otimes \Lambda^q \bar T^*.$$
The $\bar\partial_J$ complex is then the direct sum over $p$ of the Dolbeault complexes for polyvector fields:
$$\rightarrow\Omega^{0,q}(\Lambda^pT)\stackrel{\bar\partial}\rightarrow \Omega^{0,q+1}(\Lambda^pT).$$

\noindent 2. Now consider a twisted version of this example. The $\bar\partial_J$ operator is now a twisted version of the Dolbeault complex, similar to the twisted exterior derivative in Section \ref{twist}. It is easiest to describe if we choose an isotropic splitting. Then we have to add onto the Courant bracket  in equation (\ref{da}) the $H$-term, evaluated on  $(0,1)$-vectors. It gives us 
$$\bar\partial_J=\bar\partial-H^{1,2}$$
where $\bar\partial$ is the  usual Dolbeault  operator 
$$\bar\partial: \Omega^{0,q}(\Lambda^pT)\rightarrow \Omega^{0,q+1}(\Lambda^pT)$$
 and the 3--form $H^{1,2}$  acts by contraction in the $(1,0)$ factor and exterior product on the $(0,2)$ part:
 $$H^{1,2}:\Omega^{0,q}(\Lambda^pT)\rightarrow \Omega^{0,q+2}(\Lambda^{p-1}T).$$
Note that the total degree is unchanged $(q+1)+p=(q+2)+(p-1)$. 
\end{exs}

\section{Generalized K\"ahler manifolds}
\subsection{Bihermitian metrics}
Given that  a complex structure $I$ and a symplectic structure $\omega$ are both types of generalized complex structure, it makes sense to consider the analogue of a K\"ahler structure.  The compatibility of $I$ and $\omega$ is that $\omega$ should be of type $(1,1)$ with respect to $I$ (an algebraic condition) and also that  the resulting Hermitian inner product should be positive definite. The algebraic condition when translated into a condition on two generalized complex structures $J_1,J_2$ is that they should commute. We make this then a definition, and see what more we can find beyond ordinary K\"ahler metrics:
\begin{definition} A generalized K\"ahler structure on a manifold consists of a pair of commuting generalized complex structures $J_1,J_2$ such that the inner product $(J_1J_2u,v)$ is positive definite.
\end{definition}
Note that $(J_1J_2u,v)=-(J_2u,J_1v)=(u, J_2J_1v)=(u, J_1J_2v)$ so that $(J_1J_2u,v)$ does define a symmetric bilinear form.

This is a natural definition in generalized geometry, but it gives in general the notion of a {\it bihermitian metric}. This is the theorem of  Gualtieri:
\begin{thm} \label{Gthm} A generalized K\"ahler structure on a manifold gives rise to the following:
\begin{itemize}
\item
a Riemannian metric $g$
\item
two integrable complex structures $I_+,I_-$, Hermitian with respect to $g$
\item
affine connections $\nabla_{\pm}$ with skew torsion $\pm H$ which preserve the metric and the complex structure $I_{\pm}$.
\end{itemize}
Conversely, given this data, we can define a generalized K\"ahler structure, unique up  to the action of a B-field.  
\end{thm}
The proof below is based on that of \cite{MG2}. The surprising thing about this theorem is that it reveals a geometry considered long ago by  the physicists Gates, Hull and Ro\v cek \cite{GHR}.
\begin{prf}

\noindent 1.  Since the two generalized complex structures are orthogonal transformations and commute, $(J_1J_2)^2=(-1)(-1)=1$ and we have $\pm 1$ orthogonal eigenspaces $V$ and $V^{\perp}$. Splitting $u\in T\oplus T^*$ into components $u^+ + u^-$, the positive definiteness of $(J_1J_2u,v)$ means that $(u^+,u^+)-(u^-,u^-)$ is positive definite. So $V$ is positive definite and $V^{\perp}$ is negative definite. But the signature of the inner product is $(n,n)$ so $\dim V=\dim V^{\perp}=n$. This is a generalized metric as defined in Section {\ref{skewsection}} and so  already gives us  metric connections with skew torsion. 

\noindent 2. Next we find the complex structures. Since $J_1$ commutes with $J_1J_2$, it preserves the eigenspaces and defines an almost complex structure $I_+$ on $M$ by $J_1X^+=(I_+X)^+$ (on $V$, $J_2=-J_1$ and so $J_2$ defines the opposite structure). On $V^{\perp}$ $J_1=J_2$ defines similarly a complex structure $I_-$. Complexifying $T\oplus T^*$ we have a decomposition into four equidimensional subbundles corresponding to the pairs of eigenvalues of $\pm i$ of $J_1$ and $J_2$: $V^{++},V^{+-},V^{-+},V^{--}$.

So $V$, where $J_1J_2=1$ is equal to $V^{+-}\oplus V^{-+}$ and  $V^{+-}$ projects to the $+i$ eigenspace of $I_+$ on $T\otimes \C$. These subbundles are intersections of eigenspaces of $J_1,J_2$ and hence their sections are closed under Courant bracket. But the vector field part of the Courant bracket is just the Lie bracket, hence the 
$+i$ eigenspace of $I_+$ is closed under Lie bracket and hence $I_+$ is integrable. 

\noindent 3.  We need to show that the connection $\nabla$ on $V$ preserves the complex structure $I_+$, or in other words preserves the subbundle $V^{+-}$. Since this is maximal isotropic in $V\otimes \C$, we need to show that for any sections $u$ and $v$ of  $V^{+-}$  the inner product  $(\nabla_Xu,v)=0$. By the Courant bracket definition of the connection we need $([X^-,u],v)=0$. 

Decompose $X^-=x_1+x_2$ in  $V^{\perp}=V^{++}\oplus V^{--}$. Now $x_1,u$ are both in the $+i$ eigenspace of $J_1$, hence so is the Courant bracket $[x_1,u]$. Since $v$ is also in this eigenspace, which is isotropic, we have $([x_1,u],v)=0$. Similarly $x_2$ is in the $-i$ eigenspace of $J_2$, as are $u$ and $v$, so $([x_2,u],v)=0$. It follows that   $([X^-,u],v)=([x_1,u],v)+([x_2,u],v)=0$.

We now have a connection $\nabla^+$ with skew torsion which preserves the metric and the complex structure, and so is the Bismut connection.  
\vskip .25cm
\noindent 4. For the converse, assume we have the bihermitian data. Then the graph of the metric defines $V\subset T\oplus T^*$ and we use the $H$-twisted Courant bracket which defines the connection. A closed B-field will transform this to another $V$, but defining the same connections on $T$, so there is an ambiguity at this stage.

The complex structure $I_+$ splits the complexification of $V$ into $(1,0)$ and $(0,1)$ parts $V^{+-}\oplus V^{-+}$
and similarly $I_-$ splits $V^{\perp}$ complexified into $V^{++}\oplus V^{--}$. We shall prove firstly that sections of these subbundles $V^{+-}$ etc. are  closed under the Courant bracket, and then that sections of $V^{++}\oplus V^{+-}$ are also closed. Defining $J_1$ as having $+i$ eigenspace $V^{++}\oplus V^{+-}$, the closure condition will then make $J_1$ into a generalized complex structure.

So consider first $V^{+-}$. Choose local holomorphic coordinates $z_1,\dots,z_m$ with respect to $I_+$. Then elements of $V^{+-}$ can be written as 
$$\frac{\partial}{\partial z_i}+g_{i\bar j}d\bar z_j=\frac{\partial}{\partial z_i}+i\omega_{i\bar j}d\bar z_j$$
where $\omega$ is the Hermitian form.  As in the calculation in Section \ref{skewsection}, the Courant bracket
$$\left[\frac{\partial}{\partial z_i}+i\omega_{i\bar k}d\bar z_k, \frac{\partial}{\partial z_j}+i\omega_{j\bar \ell}d\bar z_{\ell}\right]=
{i(\partial \omega)}_{ji\bar\ell}d\bar z_{\ell}-H_{ji\bar\ell}d\bar z_{\ell}-H_{ji \ell}dz_{\ell}.$$
But the connection is compatible with $g$ and $I_+$ and is thus the Bismut connection, for which $H=d^c\omega$. This is a form of type $(2,1)+(1,2)$ and so has no $(3,0)$ component so $H_{ji\ell}=0$; the $(2,1)$ component is $\partial\omega$ and so the Courant bracket  vanishes.
The other three cases are similar. 

\noindent  5. Now consider $V^{++}\oplus V^{+-}$. The Courant bracket preserves sections of each component so we only have to check that $[u,v]\in V^{++}\oplus V^{+-}$ for $u$ a section of  $V^{++}$ and $v$ of $V^{+-}$. Because this is maximal isotropic in $(T\oplus T^*)\otimes \C$, we  need, as above,  the vanishing of $([u,v],w)$ for $w=w^++w^-$ a section of $V^{++}\oplus V^{+-}$.  

We just showed that the Courant bracket preserves $V^{++}$ so $[u,w^+]$ is a section of $V^{++}$ and hence $(v,[u,w^+])=0$ by isotropy. Property (\ref{cour2}) of the Courant bracket is, for any $u,v,w$,  
 $$X(v,w)=([u,v]+d(u,v),w)+(v,[u,w]+d(u,w)).$$
Since $V^{++}\oplus V^{+-}$ is isotropic this means that for our $u,v,w$ 
$$([u,v],w)+(v,[u,w])=0.$$
In particular, since $(v,[u,w^+])=0$, $([u,v],w^+)=0$. 

Similarly $w^-$ and $v$ are sections of $V^{+-}$, hence $[v,w^-]$ is a section of $V^{+-}$, so that $(u,[v,w^-])=0$ and then $([u,v],w^-)=0$.

Hence  $([u,v],w)= ([u,v],w^+)+([u,v],w^-)=0$ and $J_1$ satisfies the integrability condition; the same  argument works for $J_2$.
\end{prf}

Giving examples of generalized K\"ahler manifolds is not so straightforward as in previous notions. There are some ad hoc explicit constructions, but the existence theorem of R. Goto \cite{Go} is the most powerful method to date, and we describe this next.

\subsection{Goto's deformation theorem}
Recall that a holomorphic Poisson structure $\sigma$ gives a generalized complex structure with $E^{1,0}$ spanned by 
$$\frac{\partial}{\partial \bar z_1}, \frac{\partial}{\partial  \bar z_2},\dots,  dz_1-\sigma(dz_1), dz_2-\sigma(dz_1),\dots$$
But for $t\in \C$, $t\sigma$ is still a Poisson structure and as $t\rightarrow 0$  we get a smooth family $J_1(t)$ of generalized complex structures where $J_1(0)$ is just the standard complex structure with $E^{1,0}$
defined by 
 $$\frac{\partial}{\partial \bar z_1}, \frac{\partial}{\partial  \bar z_2},\dots,  dz_1, dz_2,\dots.$$
Goto's idea is to start with a holomorphic Poisson manifold together with a K\"ahler structure. The complex structure gives the generalized complex structure $J_1=J_1(0)$ and the symplectic structure of the K\"ahler form gives another, $J_2$. One then attempts to find, for small enough $t$, a family of generalized complex structures $J_2(t)$ which commute with $J_1(t)$ and for which $J_2(0)=J_2$.

Of course to use this to construct  generalized K\"ahler manifolds one needs to start with a compact K\"ahler holomorphic Poisson manifold, but there are a number of examples:

\begin{exs}

\noindent 1. Any algebraic surface with a holomorphic section of $\Lambda^2T$ -- this is the anticanonical line bundle. In dimension two, the integrability of $\sigma$ is automatic. Examples are $\CP^2$ blown up at  $k$ points on a cubic curve. 

\noindent 2. The Hilbert scheme of $n$ points on a Poisson surface.

\noindent 3. Any threefold with a square root $K^{1/2}$ of the canonical bundle such that  $K^{-1/2}$ has at least two sections: for example a Fano threefold of index 2 or 4. In this case take two sections $s_1,s_2$ of $K^{-1/2}$, then the Wronskian $s_1ds_2-s_2ds_1$ is a section of $T^*(K^{-1})\cong \Lambda^2T$  which satisfies the integrability condition. 

\end{exs} 
 
\begin{rmk}  In fact, any bihermitian manifold defines a holomorphic Poisson structure -- the tensor $g([I_+,I_-]X,Y)$ is the real part of a holomorphic Poisson structure with respect to either $I_+$ or $I_-$ \cite{NJH3}. It means that Poisson geometry is a central theme in this area but it would be confusing at this point to discuss the relationship between all three of these  holomorphic Poisson structures. 
\end{rmk}
 
  Here then is Goto's theorem: 
 
\begin{thm} Let $M$ be a compact K\"ahler manifold with holomorphic Poisson structure $\sigma$. Let $J_1(t)$ be the generalized complex structure defined by $t\sigma$, then for sufficiently small $t$ there exists an analytic family of  generalized K\"ahler structures $(J_1(t),J_2(t))$.
\end{thm}

\begin{prf}

The proof uses the fact that, from a linear algebra point of view, $J_1,J_2$ reduce the structure group of $T\oplus T^*$ to $U(m)\times U(m)$ -- the unitary structures on $V$ and $V^{\perp}$. In particular any two such structures are pointwise equivalent under the action of $SO(2m,2m)$. The idea is then to seek  a formal  power series $z(t)=tz_1+t^2z_2+\dots $ of skew adjoint endomorphisms of $T\oplus T^*$ such that 
\begin{itemize}
\item
 $\exp z(t)$ transforms $J_1(0)$ to $J_1(t)$ and
\item
$d(\exp z(t)e^{i\omega})=0$.
\end{itemize}
and then to 
\begin{itemize}
\item
solve the equations term-by-term and then 
\item
prove convergence by using Green's functions and harmonic  theory.
\end{itemize}
Recall from Section \ref{gcs} that $e^{i\omega}$ is the pure spinor which defines the symplectic generalized complex structure. Then  $\psi=\exp z(t)e^{i\omega}$ is again pure and satisfies $\langle\psi,\bar\psi\rangle=\langle e^{i\omega},e^{-i\omega}\rangle$ and therefore defines a $J_2(t)$. And if $d\psi=0$ it follows from Proposition \ref{int} that $J_2(t)$ satisfies the integrability condition. The last part of the process, proving convergence, is quite standard in this type of deformation theory, so we shall focus on the first part, which uses a number of features of generalized geometry. 
\vskip .25cm
\noindent 1. Let $\sigma$ be the Poisson tensor and put $a=\sigma+\bar\sigma$. This is a  section of the real exterior power $\Lambda^2T$, which is part of the bundle of skew-adjoint endomorphisms of $T\oplus T^*$, i.e. $\Lambda^2T \oplus  \End T \oplus \Lambda^2T^*$.

Let $b$ be a skew-adjoint endomorphism which preserves the generalized complex structure $J_1(0)$, which is the ordinary complex structure. Then $b$ is a real section of 
$\Lambda^{1,1}T^{1,0} \oplus \End_{\C} T\oplus \Lambda^{1,1}(T^{1,0})^*.$

For a power series $b(t)$, the composition $\exp at \circ \exp b(t)$ can (for small enough $t$) be written as $\exp z$ for a power series $z=z(t)$ and by construction its action transforms $J_1(0)$ to $J_1(t)$.  

\noindent 2. The Clifford algebra of a vector space $W$ has a filtration according to the product of generators in $W$: $\Cliff^0\subset \Cliff^2\subset \Cliff^4\subset...$ and $\Cliff^1\subset \Cliff^3\subset \Cliff^5\subset...$ (the parity is preserved because $x\cdot y+y\cdot x=2(x,y)1$ is of degree zero). We saw in Section \ref{spin} that $\{a\in \Cliff(W): [a,W]\subseteq W\, {\mathrm {and}} \,\,  a=-a^t\}$ is isomorphic under the spin representation to the skew-adjoint endomorphisms of $W$, so we consider $z$ as lying in $\Cliff^2$ and exponentiation in the group is exponentiation in the Clifford algebra. It follows that in $\Cliff(T\oplus T^*)$ we have 
$$e^{-z}\Cliff^1e^z \subseteq \Cliff^1.$$

\noindent 3. We need to solve $d(e^ z \cdot \psi)=0$  where $\psi=e^{i\omega}$ so we consider the operator $e^{-z} d\, e^z$ for $z$ a section of $\Cliff^2$ (note that if $z=B$  in $\Lambda^2T^*\subset \Lambda^2T \oplus  \End T \oplus \Lambda^2T^*$, then $e^{-z} d\, e^z$ is the twisted differential $d+H$ where $H=dB$ which we met in Section  \ref{twist}, but here we need the general case). 

We use the formula
$$d\varphi=\sum_idx_i\wedge{\mathcal L}_{\partial/\partial x_i}\varphi.$$
(This is trivially true on functions and the right hand side has the obvious property $d(\alpha\wedge\beta)=d\alpha\wedge \beta+(-1)^p\alpha\wedge d\beta$.)
Then
\begin{eqnarray*}
e^{-z} d\, e^z \cdot \varphi &=&e^{-z}  \sum_idx_i\wedge{\mathcal L}_{\partial/\partial x_i}( e^z \cdot \varphi)\\
&=&\sum_i(e^{-z}\cdot dx_i\cdot  e^z)\cdot e^{-z}{\mathcal L}_{\partial/\partial x_i} (e^z \cdot \varphi)
\end{eqnarray*}
where we have rewritten the exterior product by a 1-form as a Clifford product. This expands to 
$$\sum_i(e^{-z}\cdot dx_i\cdot  e^z)((e^{-z}{\mathcal L}_{\partial/\partial x_i}e^z) \varphi+{\mathcal L}_{\partial/\partial x_i}  \varphi).$$
Now $u_i=e^{-z}\cdot dx_i\cdot  e^z$ is a local section of $\Cliff^1$ and $a_i=(e^{-z}{\mathcal L}_{\partial/\partial x_i}e^z)$ is a section of the Lie algebra bundle in $\Cliff^2$, so we have
\begin{equation}
e^{-z} d \,e^z \varphi =\sum_i u_i\cdot a_i\cdot\varphi+u_i\cdot {\mathcal L}_{\partial/\partial x_i}  \varphi
\label{dtwist}
\end{equation}

\noindent 4. We need to consider the action of $\exp z$ on the pair $J_1(0),J_2(0)$ and it is convenient to see this via the pair of pure spinors $(\psi,\varphi)=(e^{i\omega}, dz_1\wedge\dots\wedge dz_m)$ -- the first is global, the second only local. The Lie derivative action on $T\oplus T^*$ preserves the inner product and the pair $(\psi,\varphi)$ at each point lie in an $SO(2m,2m)$ orbit, so
$${\mathcal L}_{\partial/\partial x_i} (\psi,\varphi)=(c_i\cdot\psi,c_i\cdot \varphi)$$
for some local section $c_i$ of $\Cliff^2$. From equation (\ref{dtwist}) we then have 
$$e^{-z} d \,e^z (\psi,\varphi)=(h\cdot\psi,h\cdot\varphi)$$
for some $h=\sum_i u_i\cdot(a_i+c_i)$ a section of $\Cliff^3$. 

But by our choice of $z$,  $e^z \varphi$ defines a generalized complex structure so by the integrability criterion in Proposition \ref{int} we have $d(e^z \varphi)=w\cdot e^ z \cdot \varphi$ for some $w$ a local section of $\Cliff^1=T\oplus T^*$. Putting $v=e^{-z}\cdot w \cdot e^z$ this gives us two equations: an algebraic condition on $h$, a section of $\Cliff^3$, that $v\cdot \varphi=h\cdot \varphi$ for some $v$ a section of $\Cliff^1$, and the differential equation
\begin{equation}
e^{-z} d( e^z  \psi)=h\cdot\psi.
\label{deq}
\end{equation}
The important thing to note is that  there is an $h$ satisfying these conditions for any $b(t)=tb_1+t^2b_2+\dots...$.

\noindent 5.  We now need to identify the objects on the right hand side of this equation, which is the second item in the following lemma:

 \begin{lemma} \label{hlemma} Let $b\in\Cliff^2$ and $h\in \Cliff^3$ be real elements, then

\noindent (i) if  $b$ preserves $dz_1\wedge\dots\wedge dz_m$ up to a scalar multiple, then $b\cdot e^{i\omega}\in (\Lambda^0+\Lambda^{1,1})\wedge e^{i\omega}$.

\noindent (ii) If  $h$ satisfies $h\cdot dz_1\wedge\dots\wedge dz_m=v\cdot dz_1\wedge\dots\wedge dz_m$ for some $v\in \Cliff^1$ then $h\cdot e^{i\omega}\in (\Lambda^1+\Lambda^{2,1}+\Lambda^{1,2})\wedge e^{i\omega}$.
\end{lemma}

\begin{lemprf}

\noindent (i) If $X$ is any tangent vector, $i_Xe^{i\omega}=i(i_X\omega)\wedge e^{i\omega}$ so the action of $b$ can always be realized as the  exterior product by some 0-form plus 2-form. The $(1,1)$ forms annihilate $dz_1\wedge\dots\wedge dz_m$ and so $\Lambda^{1,1}\wedge e^{i\omega}$ is in the image. Terms in $(\Lambda^{2,0}+\Lambda^{0,2})\wedge e^{i\omega}$  can arise from real  linear combinations of the real and imaginary parts of
$$\frac{\partial}{\partial z_i}\cdot \frac{\partial}{\partial z_j}\qquad\frac{\partial}{\partial z_i}\cdot d\bar z_j \qquad  d\bar z_i\cdot
d\bar z_j$$
But applied to $dz_1\wedge\dots\wedge dz_m$ these give non-zero terms with the  respective $(p,q)$ types $(m-2,0), (m-1,1), (m,2)$  and so do not preserve the pure spinor.

\noindent (ii) If $v\in \Cliff^1$ then $v\cdot dz_1\wedge\dots\wedge dz_m$ is of type $(m-1,0)+(m,1)$. As in (i), the action of $h$ on $e^{i\omega}$ is the exterior product of   some 1-form plus 3-form. Forms of type $(2,1)+(1,2)$ annihilate $dz_1\wedge\dots\wedge dz_m$ by exterior multiplication so $(\Lambda^{2,1}+\Lambda^{1,2})\wedge e^{i\omega}$ is in the image. This time terms in $(\Lambda^{3,0}+\Lambda^{0,3})\wedge e^{i\omega}$ can arise from  combinations  of
$$\frac{\partial}{\partial z_i}\cdot \frac{\partial}{\partial z_j}\cdot \frac{\partial}{\partial z_k}\qquad\frac{\partial}{\partial z_i}\cdot \frac{\partial}{\partial z_j}\cdot d\bar z_j \qquad \frac{\partial}{\partial z_i}\cdot d\bar z_j \cdot d\bar z_k\qquad  d\bar z_i\cdot
d\bar z_j\cdot d\bar z_k$$
but these applied to $dz_1\wedge\dots\wedge dz_m$ give non-zero forms with   the $(p,q)$ types $(m-3,0)$, $(m-2,1),(m-1,2),(m,3)$ and not  the required types $(m-1,0),(m,1)$.

It is easy to see that when $b$ or $h$ are real, we can realize the action by a  real form. 
\end{lemprf}

\noindent 6. We now begin the term-by-term solution to the equation $d( e^z \psi)=0$. Recall that 
$$e^{z(t)}=e^{ta}e^{b(t)}=1+t(a+b_1)+\dots$$ and so
$$e^{-z} d( e^z  \psi)=td((a+b_1)\psi)+o(t^2)$$
so the first task is to solve $d((a+b_1)\psi)=0$ for $b_1$. If we put all the $b_i=0$ in Equation (\ref{deq}), and use item (ii) in Lemma \ref{hlemma} then we see that 
$$d(a\psi)\in (\Omega^1+\Omega^{2,1}+\Omega^{1,2})\wedge \psi.$$
 But exterior product with $\psi=e^{i\omega}$ is invertible  and commutes with $d$, so we have an exact form in $\Omega^1+\Omega^{2,1}+\Omega^{1,2}$ and, to find $b_1$, from item (i) in the  lemma we want  this  to lie in $d(\Omega^0+\Omega^{1,1})$. 
 
The 1-form part is obvious. Write the   component in $\Omega^{2,1}+\Omega^{1,2}$ as the exact form $d(\alpha^{2,0}+\alpha^{1,1}+\alpha^{0,2})$ then $\bar\partial \alpha^{0,2}=0$ so $\partial \alpha^{0,2}$ is $\bar\partial$-closed and $\partial$-exact. So by the $\partial\bar\partial$-lemma (valid for K\"ahler manifolds) we can write 
$$\partial \alpha^{0,2}=\partial \bar\partial\theta^{0,1}$$
and so the $(1,2)$ component of the exact form can be written as 
$$\bar\partial \alpha^{1,1}+\partial \alpha^{0,2}=\bar\partial(\alpha^{1,1}- \partial\theta^{0,1})$$
and then $\alpha^{1,1}- \partial\theta^{0,1}-\overline{\partial\theta^{0,1}}$ is the required $(1,1)$ form. 

\noindent 7. In general suppose we have inductively found $b_1,\dots,b_{k-1}$ so that $d(e^{z(t)}\psi)_i=0$ for $i<k$, then
$$(e^{-z(t)} d( e^{z(t)}  \psi))_k=\sum_{i+j=k}(e^{-z(t)})_jd( e^{z(t)} \psi)_i=(d(e^{ z(t)} \psi))_k.$$
Then using (\ref{deq}) and the  $\partial\bar\partial$-lemma again we can solve for $b_k$.
 \end{prf}
 
 \subsection{Deformation of the  bihermitian structure}\label{deform}
 
According to  Gualtieri's theorem, Goto's result gives a pair of complex structures $I_+,I_-$ starting from a K\"ahlerian holomorphic Poisson manifold. One may ask where these are coming from, and we can answer this at the infinitesimal level. 

Recall that given a smooth family $I(t)$ of complex structures  the derivative $I'$ satisfies  $ I' I+I  I'=0$ and defines a $\bar\partial$-closed section of $T^{1,0}\otimes (T^*)^{0,1}$ which gives the Kodaira-Spencer class in the Dolbeault cohomology group $[I']\in H^1(M,T)$.  More concretely, if we write 
$$\left(\frac{\partial}{\partial \bar z_j}\right)'= \alpha_{\bar j k}\frac{\partial}{\partial  z_k}+\beta_{\bar j \bar k}\frac{\partial}{\partial  \bar z_k}$$
then the Kodaira-Spencer class is represented by 
$$ \alpha_{\bar j k}\frac{\partial}{\partial  z_k}\,d\bar z_j$$
We defined the complex structure $I_+$ in Theorem \ref{Gthm}  in terms of the intersection of the $+i$ eigenspace of $J_1$ and the $-i$ eigenspace of $J_2$: in the K\"ahler case this  gives sections of $(T\oplus T^*)\otimes \C$ spanned by terms of the form 
$$u_j=\frac{\partial}{\partial \bar z_j}+i\omega_{\bar j k}dz_k.$$
Deform this and we are looking for sections $u_j(t)$ such that 
$u_j\cdot \psi=0$ and $u_j\cdot \varphi=0$. So differentiating with respect to $t$ at $t=0$,
$u'_j\cdot \psi+u_j\cdot \psi'=0$ and $ u_j' \cdot \varphi+u_j\cdot \varphi'=0$. 
Now $\varphi(t)=e^{at}dz_1\wedge\dots\wedge dz_m$, so 
$$\varphi'=a\cdot dz_1\wedge\dots\wedge dz_m= \sigma^{pq}i_{\partial/\partial z_p}i_{\partial/\partial z_q}dz_1\wedge\dots\wedge dz_m.$$
and hence, considering the $(m-1,0)$ component of the equation $  u_j' \cdot \varphi+u_j\cdot \varphi'=0$ we get
 $$\alpha_{\bar j k}i_{\partial/\partial z_k}dz_1\wedge\dots\wedge dz_m+i\omega_{\bar j k}dz_k\wedge \sigma^{p q}i_{\partial/\partial z_{p}}i_{\partial/\partial z_q}dz_1\wedge\dots\wedge dz_m=0.$$
 But
 $$i_{\partial/\partial z_{p}}(d z_k\wedge i_{\partial/\partial z_{q}}dz_1\wedge\dots\wedge dz_m)=\delta_{pk} i_{\partial/\partial z_{q}}dz_1\wedge\dots\wedge dz_m-dz_k\wedge i_{\partial/\partial z_{p}}i_{\partial/\partial z_q}dz_1\wedge\dots\wedge dz_m$$
and 
$$d z_k\wedge i_{\partial/\partial z_{q}}dz_1\wedge\dots\wedge dz_m=\delta_{qk}dz_1\wedge\dots\wedge dz_m.$$
It follows that 
$$\alpha_{\bar j k}=-2i\omega_{\bar j \ell}\sigma^{\ell k}.$$
More invariantly, the K\"ahler form $\omega$  defines a Dolbeault class in $H^1(M,T^*)=H^{1,1}$ and the Poisson tensor $\sigma$ lies in $H^0(M,\Lambda^2T)$; then the Kodaira-Spencer class in $H^1(M,T)$ for the deformation $I_+(t)$ is the cup product combined with contraction ${\sigma}\omega$. If we do the same for $I_-(t)$ we get the negative of this class.

\begin{exs}

\noindent 1. For $\CP^2$, the bundle $\Lambda^2 T$ is isomorphic to ${\mathcal O}(3)$, and a holomorphic Poisson structure is defined by a section of this, which vanishes on a cubic curve. As mentioned earlier, blowing up $r$ points on this curve, the proper transform is again an anticanonical divisor and so defines a Poisson structure. If $E_1,\dots, E_{r}$ are the cohomology classes of the exceptional curves on the blow-up and $H$ the pull-back of a  generator of $H^2(\CP^2,\Z)$ then a K\"ahler form $\omega$ has a cohomology class of the form
$$[\omega]=aH+\sum_1^rc_iE_i$$
where $a>0$ and $c_i<0$. In this case the Kodaira-Spencer class for  $I_+$ above consists of moving the points along the cubic curve with velocities (relative to a trivialization of the tangent bundle of the elliptic curve) proportional to $c_i$. For $I_-$, since $c_i$ have the same sign, the points all move in the opposite direction.

\noindent 2. If the cohomology class of the K\"ahler form is a multiple of $c_1(T)$  then, as we shall see in the next lecture, the Kodaira-Spencer class is zero.  This is consistent with some concrete constructions of one-parameter families of bihermitian metrics on Fano manifolds in \cite{NJH4} and \cite{MG3}, where $I_+(t)$ and $I_-(t)$ are all equivalent under a diffeomorphism.

\end{exs}

One corollary of Goto's theorem is that a Kodaira-Spencer class of the form just described is {\it unobstructed} -- there exists a one parameter family integrating it. This is reminiscent of the Tian-Todorov result on Calabi-Yau manifolds. (In fact the analogy is close since Goto's theorem is a deformation theorem of generalized Calabi-Yau manifolds: the generalized complex structure $J_2(t)$ is defined by a closed pure spinor.) One can see this, however, directly without the language of generalized geometry:

\begin{prp} Let $M$ be a compact complex  manifold with K\"ahler form $\omega$ and holomorphic Poisson tensor $\sigma$. Then the Kodaira-Spencer class $\sigma\omega\in H^1(M,T)$ is unobstructed.
\end{prp}

\begin{prf} For the deformation of a complex structure, one looks for  a global $\phi\in \Omega^{0,1}(T)$ which satisfies the equation
$$\bar\partial\phi =[\phi,\phi]$$
where the bracket is the Lie bracket on  vector fields together with  exterior product on $(0,1)$-forms, and, as above, one solves this (if possible) term-by-term for a series $\phi(t)=t\phi_1+t^2\phi_2+\dots$ where $\phi_1$ represents the initial Kodaira-Spencer class. 

The coefficient of $t^2$ requires a solution for $\phi_2$ of 

\begin{equation}
\bar\partial \phi_2=[\phi_1,\phi_1]. 
\label{ob}
\end{equation}
 Now locally we have the K\"ahler potential 
$$\omega_{\bar j \ell}=\frac{\partial^2 f}{\partial  z_{\ell} \partial \bar z_j}$$
and then writing $f_{\bar j}=\partial f/\partial \bar z_j$,
$$\phi_1=\omega_{\bar j \ell}\sigma^{\ell k}\frac{\partial}{\partial z_k}\,d\bar z_j=\sigma^{\ell k}\frac{\partial^2 f}{\partial  z_{\ell} \partial \bar z_j}\frac{\partial}{\partial z_k}\,d\bar z_j=\sigma(\partial f_{\bar k})d\bar z_k.$$
Now $\sigma$ is holomorphic and so is insensitive to $\bar\partial$ operators, therefore 
$$[\phi_1,\phi_1]=\sigma(\partial\{f_{\bar k},f_{\bar \ell}\})d\bar z_k\wedge d\bar z_{\ell}$$
using the integrability of $\sigma$. 

The Poisson bracket expression $\{f_{\bar k},f_{\bar \ell}\}d\bar z_k\wedge d\bar z_{\ell}$ looks local  but it is 
$$\sigma^{ij}f_{i\bar k}f_{j \bar \ell}=\sigma^{ij}\omega_{i\bar k}\omega_{j\bar \ell}$$
or more invariantly $i_{\sigma}(\omega\wedge\omega)$. Thus $\partial i_{\sigma}(\omega\wedge\omega)=[\phi_1,\phi_1]$ is $\bar\partial$-closed and $\partial$-exact and so, by the $\partial\bar\partial$-lemma 
$$\partial i_{\sigma}(\omega\wedge\omega)=\bar\partial\partial\alpha$$
for some $(0,1)$-form $\alpha$. It follows that 
$$[\phi_1,\phi_1]=\sigma(\partial  i_{\sigma}(\omega\wedge\omega))=\sigma(\bar\partial\partial\alpha)=\bar \partial( \sigma(\partial\alpha))$$
and we take $\phi_2=\sigma(\partial\alpha)$ to solve Equation (\ref{ob}).
\vskip .25cm
The coefficient of $t^3$ typifies the inductive process: we need 
$$\bar\partial \phi_3=[\phi_1,\phi_2]+[\phi_2,\phi_1].$$
Write $\alpha=g_{\bar \ell}d\bar z_{\ell}$. then $[\phi_1,\phi_2]=\sigma(\partial \{f_{\bar k},g_{\bar \ell}\})d\bar z_k\wedge d\bar z_{\ell}$ and 
$$\{f_{\bar k},g_{\bar \ell}\}d\bar z_k\wedge d\bar z_{\ell}=i_{\sigma}(\omega\wedge \partial\alpha)$$
and we proceed as above.
\end{prf}

\section{Generalized holomorphic bundles}
\subsection{Basic features}
The analytic viewpoint of a holomorphic vector bundle $V$ on a complex manifold was established by Malgrange (see \cite{DK} Section 2.2.2 for a simple proof).  The existence of  sufficiently many local holomorphic sections is equivalent to the existence of a differential operator $\bar\partial_A:\Omega^0(V)\rightarrow \Omega^{0,1}(V)$ such that $\bar\partial_A(fs)=\bar\partial f s+f\bar\partial_A s$ and such that the standard extension to forms $\bar\partial_A^2:\Omega^0(V)\rightarrow \Omega^{0,2}(V)$ vanishes. Gualtieri introduced the analogous concept in generalized complex geometry:
\begin{definition} A generalized holomorphic bundle on a generalized complex manifold $(M,J)$ is a vector bundle $V$ with a differential operator $\bar D:C^{\infty}(V)\rightarrow C^{\infty}(V\otimes E^{0,1})$ such that for a smooth function $f$ and section $s$ 
\begin{itemize}
\item
$\bar D(fs)=\bar\partial_J \!f s+f\bar Ds$
\item
$\bar D^2: C^{\infty}(V)\rightarrow C^{\infty}(V\otimes \Lambda^2E^{0,1})$ vanishes.
\end{itemize}
\end{definition}

Given a local trivialization $s_1,\dots,s_k$ of $V$ we obtain a ``connection matrix" $A_{ij}$ with values in $E^{0,1}$ defined by
$$\bar D s_i=A_{ji}s_j$$
and then the condition $\bar D^2=0$ is $\bar\partial_J A+A\cdot A=0$.

\begin{rmk} For an ordinary holomorphic bundle we have a Dolbeault complex 
$$\rightarrow \Omega^{0,p}(V)\stackrel{\bar\partial}\rightarrow\Omega^{0,p+1}(V)\rightarrow $$
and by the same token there is a generalized version
$$\rightarrow C^{\infty}(V\otimes \Lambda^p E^{0,1})\stackrel{\bar D}\rightarrow C^{\infty}(V\otimes \Lambda^{p+1} E^{0,1})\rightarrow .$$
\end{rmk}

 A universal example of a generalized holomorphic bundle  is the canonical bundle of the generalized complex structure -- the 
subbundle $K\subset \Lambda^*T^*\otimes\C$ of multiples of pure spinors whose annihilator is $E^{1,0}$. To see this recall Proposition \ref{int}, where we saw that $J$ was integrable if and only if for any local non-vanishing section $\psi$ of $K$, we have $d\psi=w\cdot\psi$ for some local section $w$ of $(T\oplus T^*)\otimes \C$. The $(1,0)$ component of $w$ annihilates $\psi$ and the $(0,1)$ component  is unique because $E^{1,0}\cap E^{0,1}=0$, so we may as well assume $w$ lies in $E^{0,1}$. 

A global section $s$ of $K$ can be written locally as $f\psi$ and we define
$$\bar D s=(\bar\partial_Jf) \psi+f w\cdot\psi.$$
If $\psi_1=g\psi$ is another local section, then $w_1=g^{-1}dg+w$, $f=f_1g$ and  one can easily check that $\bar D s$  is well-defined. 

We need also the condition $\bar D^2=0$ which means we need to prove that $\bar\partial_J w=0$.  Let $u=X+\xi, v=Y+\eta$ be sections of $E^{1,0}$, then $u\cdot\psi=0=v\cdot\psi$ and 
$${\bf L}_u\psi=d(u\cdot\psi)+u\cdot d \psi=u\cdot w\cdot\psi=2(u,w)\psi$$
and
$${\bf L}_v{\bf L}_u\psi=(2Y(u,w)+4(u,w)(v,w))\psi.$$
Hence
$$({\bf L}_u{\bf L}_v-{\bf L}_v{\bf L}_u)\psi=2(X(v,w)-Y(u,w))\psi.$$
But $({\bf L}_u{\bf L}_v-{\bf L}_v{\bf L}_u)\psi={\bf L}_{[u,v]}\psi$ (see the proof of Proposition \ref{jac}) and since, by integrability of $J$, $[u,v]\cdot \psi=0$ we also have 
${\bf L}_{[u,v]}\psi=2([u,v],w)\psi$.
Hence  $([u,v],w)=X(v,w)-Y(u,w)$ which from the definition (\ref{da}) is $\bar\partial_J w=0$.
 
From this we see also  that a generalized Calabi-Yau manifold, which we have defined as having  a global closed $\psi$, can also be thought of as being defined by a global non-vanishing generalized holomorphic section of the canonical bundle. 

For the specific examples of generalized complex structures --  symplectic, complex, holomorphic Poisson -- we can determine what  a generalized holomorphic bundle means:

\begin{exs}

\noindent 1. On a symplectic manifold $E^{0,1}$ is spanned by terms
$$\frac{\partial}{\partial x_j}+i\omega_{jk}dx_k$$
or equivalently, inverting $\omega_{ij}$, by
$$dx_i-i\omega_{ij}\frac{\partial}{\partial x_j}.$$
Thus the 1-form part of the connection matrix for a generalized holomorphic bundle 
$$A_i(dx_i-i\omega_{ij}\frac{\partial}{\partial x_j})$$
defines an ordinary connection and $\bar D^2=0$ implies it is a flat connection.

\noindent 2. Now consider a complex manifold considered as a generalized complex manifold. On might think that  generalized holomorphic bundles are just  ordinary holomorphic bundles, but this is not quite the full picture, although they do provide examples.  

Since  $E^{0,1}=\bar T^*\oplus T$ (where now $T$ is the holomorphic tangent bundle)
$$\bar Ds=\bar\partial_A s+\phi s=\left(\frac{\partial s}{\partial \bar z_{j}}+A_{\bar j}s\right) d\bar z_{j}+\phi^{k} s \frac{\partial }{\partial z_{k}}$$
and $\bar D^2=0$ implies 
$$\bar\partial_A^2=0 \in \End V\otimes \Lambda^2\bar T^*,\quad \bar\partial_A \phi=0 \in \End V\otimes \bar T^*\otimes T,\quad \phi^2=0\in \End V\otimes \Lambda^2T.$$
The first condition gives $V$ the structure of a holomorphic vector bundle, the second says that $\phi$ is a holomorphic section of $\End V\otimes T$, and the third  is an algebraic condition on $\phi$. Since
$$\phi^2=\frac{1}{2}[\phi^{j},\phi^{k}]\frac{\partial }{\partial  z_{j}}\wedge \frac{\partial }{\partial  z_{k}}$$
this ``integrability" condition is   $[\phi^{j},\phi^{k}]=0$.

We call these {\it co-Higgs bundles}. A Higgs bundle in the sense of C.Simpson \cite{S} is the same definition with $T$ replaced by $T^*$.

\noindent 3. The generalized complex structure determined by a holomorphic Poisson tensor has 
$E^{0,1}$ spanned  by 
$$\frac{\partial}{\partial  z_1}, \frac{\partial}{\partial   z_2},\dots,  d\bar z_1-\bar\sigma(d\bar z_1), d\bar z_2-\bar\sigma(d\bar z_1),\dots$$
and the $\bar\partial_J$ operator is 
$$\bar\partial_Jf=\bar\partial f+\sigma(\partial f)-\bar\sigma(\bar\partial f).$$
We then write $\bar D$ as 
$$\bar Ds=\left(\frac{\partial s}{\partial \bar z_{j}}+A_{\bar j}s\right) (d\bar z_{j}-\bar\sigma(d\bar z_j))+\frac{\partial s}{\partial z_{j}}\sigma(dz_j)+\phi^{k} s \frac{\partial }{\partial z_{k}}.$$
Again $A_{\bar j}$ defines a holomorphic structure on $V$ and then, in a local holomorphic basis, the operator is 
$$\bar Ds=\frac{\partial s}{\partial z_{j}}\sigma(dz_j)+\phi^{k} s \frac{\partial }{\partial z_{k}}.$$
This is a first order holomorphic differential operator from $V$ to $V\otimes T$ whose symbol is $1\otimes \sigma:V\otimes T^*\rightarrow V\otimes T$. We can define the action of a local holomorphic function $f$ on a section $s$ by
$$f\cdot s=\langle \bar D s,df\rangle$$
and then the condition $\bar D^2=0$ says that
$$g\cdot f\cdot s-f\cdot g\cdot s=\{g,f\}\cdot s.$$
Such a holomorphic bundle is called a  {\it Poisson module}.

\begin{rmk} Note that, given a co-Higgs bundle $(V,\phi)$ we can define an action of $f$ on a local section by 
$$f\cdot s=\phi(df)s$$
and then the $\phi^2=0$ condition says that $g\cdot f\cdot s-f\cdot g\cdot s=0$. We can thus interpret a co-Higgs bundle as a Poisson module for the zero Poisson structure.
\end{rmk}
\end{exs}

In the next two sections we shall examine the last two examples in more detail.

\subsection{Co-Higgs bundles}\label{coH}
A co-Higgs bundle is, as we have seen, defined by  a pair consisting of a holomorphic vector bundle $V$ and an endomorphism $\phi$, twisted by the tangent bundle. When studying such pairs one usually imposes a stability condition in order to construct a Hausdorff moduli space. On a K\"ahler manifold one can define the degree of a line bundle and the slope $\deg \Lambda^kV/k$ of a vector bundle $V$ of rank $k$. The stability condition for Higgs bundles in \cite{S} is that the slope of any $\phi$-invariant torsion-free subsheaf should be less than the slope of $V$. This makes perfectly good sense whether one takes the tangent bundle or the cotangent bundle but the manifolds which support such stable objects are quite different. In one dimension, for example,  the main interest in the case of Higgs bundles lies with  genus $g>1$, for in that case there is a link with  representations of the fundamental group. In the co-Higgs case there are no stable objects with $\phi\ne 0$ in higher genus. The point is that given a section $s$ of the $g$-dimensional space $H^0(M,K)$ of differentials, if $\phi\in H^0(M,\End V\otimes K^*)$ then $\phi s$ is an endomorphism which commutes with $\phi$ and stable objects do not have any of these other than the scalars. 

\begin{exs}

\noindent 1. In rank one, a co-Higgs bundle is just a line bundle $V=L$ together with a vector field $\phi$.

\noindent 2. If $V={\mathcal O}\oplus T$ then there is a canonical co-Higgs structure where $\phi(\lambda,X)=(X,0)$. Since the trivial bundle is invariant, we require $\deg T>0$ and $T$ itself to be stable for stability of the co-Higgs bundle.
\end{exs}

We refer the reader to the forthcoming Oxford DPhil thesis of Steven Rayan for more  results about co-Higgs bundles, but here we shall give some examples on projective spaces slightly more interesting  than those above. 

\begin{ex} Consider $M=\CP^m=\PP(W)$. The tangent bundle fits into the Euler sequence of holomorphic bundles 
$$0\rightarrow  {\mathcal O}\rightarrow         W(1)\rightarrow T\rightarrow 0$$
from which we obtain $W\cong H^0(\CP^m,T(-1))$. We also see that $\Lambda^m T\cong {\mathcal O}(m+1)$ and hence
$$T^*\cong \Lambda^{m-1}T\otimes \Lambda^mT^*\cong \Lambda^{m-1}T(-(m+1)).$$
This means that
$$T\otimes T^*(1)\cong T\otimes \Lambda^{m-1}T(-m))=T(-1)\otimes  \Lambda^{m-1}(T(-1))$$
and from the $(m+1)$-dimensional space of sections of $T(-1)$ we can construct by  tensor and exterior product many sections, not just scalars, of $T\otimes T^*(1)$. Take one, $\psi$, and a section $w$ of $T(-1)$. Then set 
$$\phi= \psi\otimes w\in H^0(\CP^m,\End T\otimes T).$$
By construction, $\phi^2= [\psi,\psi]\otimes w\wedge w=0$, and the tangent bundle itself is stable so this gives plenty of examples of co-Higgs bundles on projective space.
\end{ex}
\vskip .25cm
 The simplest concrete example, where we can write down the moduli space, is the case of the bundle $V={\mathcal O}\oplus {\mathcal O}(-1)$ on $\CP^1$. Since $\Lambda^2T=0$, there is no integrability condition on $\phi$ in one dimension.

Here the tangent bundle is ${\mathcal O}(2)$ and so we must have 
$$\phi=\pmatrix{a & b\cr
                             c & -a}$$
                             where $a,b,c$ are sections of ${\mathcal O}(2), {\mathcal O}(3), {\mathcal O}(1)$ respectively. Since the degree and rank of $V$ are coprime, there are no semi-stable bundles which means, in this one-dimensional case, that the moduli space is smooth. We first define a canonical six-dimensional complex manifold. We denote by $p:T\CP^1\rightarrow \CP^1$ the projection and $\eta \in H^0(T\CP^1,p^*T)$ the tautological section. This is  a section of $p^*{\mathcal O}(2)$. Now define
                              $${\mathcal M}=\{(x, s)\in T\CP^1\times H^0(\CP^1,{\mathcal O}(4)): \eta^2(x)=s(p(x))\}.$$
\begin{prp} ${\mathcal M}$ is naturally isomorphic to the moduli space of stable rank 2 trace zero co-Higgs bundles of  degree $-1$ on $\CP^1$.
\end{prp}     
\begin{prf}

\noindent 1. First note that any vector bundle on $\CP^1$ is a sum of line bundles by the Birkhoff-Grothendieck theorem.       If the decomposition is      ${\mathcal O}(m)\oplus {\mathcal O}(-1-m)$, then 
 $a,b,c$ are sections of ${\mathcal O}(2), {\mathcal O}(2m+3), {\mathcal O}(1-2m)$.
If $c$ is zero, then ${\mathcal O}(m)$ is $\phi$-invariant, so by stability $m< -1/2$. If $m=-1$, then by changing the order of the subbundles we are in the same situation. If $m\le -2$ then $b$ is a section of a line bundle of negative degree and so vanishes -- then the invariant subbundle ${\mathcal O}(-1-m)$ contradicts stability. Thus the vector bundle in this moduli space is always $V={\mathcal O}\oplus {\mathcal O}(-1)$.

\noindent 2. Since  $c$ is a non-zero section of ${\mathcal O}(1)$ it vanishes at a distinguished point $z=z_0$. Then $a(z_0)$ is a point $x$ in the total space of ${\mathcal O}(2)=T\CP^1$. It is well-defined because an automorphism of ${\mathcal O}\oplus {\mathcal O}(-1)$ is defined by 
$$\pmatrix{A & B\cr
                             0 & C}$$
 where $A,B,C$ are sections of ${\mathcal O}, {\mathcal O}(1), {\mathcal O}$ and where $c=0$ the action on $a$ is trivial.

\noindent 3. The determinant $\det \phi$ is a section $s$ of ${\mathcal O}(4)$, and at $z=z_0$, $c$ vanishes so we have $\det\phi(z_0)=-a(z_0)^2$, so set $p=-\det\phi$. This defines a map from the moduli space to ${\mathcal M}$.

\noindent 4. In the reverse direction, choose an affine parameter $z$ on $\CP^1$ such that the point $x$ maps to $z=0$ and write $p(z)=a_0^2+zb(z)$ where $b(z)$ is a cubic polynomial. Then $\eta^2(0)=a_0^2$ so $\eta(0)=\pm{a_0}$ and 
$$\pmatrix{\eta(0) & b(z)\cr
                             z & -\eta(0)}$$
is a representative Higgs field. 
\vskip .25cm
Note that ${\mathcal M}$ is a fibration of elliptic curves $y^2=c_0+c_1z+\dots +c_4z^4$ over the five-dimensional vector space of coefficients $c_0,\dots,c_4$. We shall see this again when we consider the B-field action in the next lecture. 
\end{prf}

\begin{rmk} We saw in Section \ref{dbarsec} that the $\bar\partial_J$ complex for an ordinary complex structure was defined by 
$$\Omega^{0,q}(\Lambda^pT)\stackrel{\bar\partial}\rightarrow \Omega^{0,q+1}(\Lambda^pT).$$
For a co-Higgs bundle $(V,\phi)$ the $\bar D$ complex is defined by $\bar\partial+\phi$ where 
$$\bar\partial: \Omega^{0,q}(V\otimes \Lambda^pT)\rightarrow \Omega^{0,q+1}(V\otimes \Lambda^pT)$$
and 
$$\phi:  \Omega^{0,q}(V\otimes \Lambda^pT)\rightarrow \Omega^{0,q}(V\otimes \Lambda^{p+1}T).$$
\end{rmk} 
Note that the total degree $p+(q+1)=(p+1)+q$ is preserved. It is easy to see that 
the cohomology of the $\bar D$ complex is the hypercohomology of the complex of sheaves 
$$\dots \rightarrow {\mathcal O}(V\otimes \Lambda^pT)\stackrel{\phi}\rightarrow {\mathcal O}(V\otimes \Lambda^{p+1}T)\rightarrow\dots$$

\subsection{Holomorphic Poisson modules}

We observed that a holomorphic Poisson module is a holomorphic vector bundle $V$ with a 
 first order holomorphic linear differential operator 
$$\bar D:{\mathcal O}(V)\rightarrow {\mathcal O}(V\otimes T)$$
whose symbol is $1\otimes\sigma: V\otimes T^*\rightarrow V\otimes T$. 
Relative to a local holomorphic basis $s_i$ of $V$, $\bar D$ is defined by a ``connection matrix" $A$ of vector fields:
$$\bar Ds_i= s_j\otimes A_{ji}.$$
When $\sigma$ is non-degenerate it identifies $T$ with $T^*$ and then $\bar D$ is a flat holomorphic connection.

\begin{ex} If $X=\sigma(df)$ is the Hamiltonian vector field of $f$ then the Lie derivative ${\mathcal L}_X$ acts on tensors  but the action in general  involves the second derivative of $f$. However for the canonical line bundle $K=\Lambda^nT^*$ we have 
$${\mathcal L}_X(dz_1\wedge\dots\wedge dz_n)=\frac{\partial X_i}{\partial z_i} (dz_1\wedge\dots\wedge dz_n)$$
and, since $\sigma^{ij}$ is skew-symmetric,
$$\frac{\partial X_i}{\partial z_i}=\frac{\partial}{\partial z_i}\left(\sigma^{ij}\frac{\partial f}{\partial z_j}\right)=\frac{\partial \sigma^{ij}}{\partial z_i}\frac{\partial f}{\partial z_j}$$
which involves only the first  derivative of $f$. Thus 
$$\{f,s\}={\mathcal L}_Xs=\langle \bar Ds, df\rangle$$
defines a first order operator. The second condition for a Poisson module follows from the integrability of the Poisson structure: since $\sigma(df)=X, \sigma(dg)=Y$ implies $\sigma(d\{f,g\})=[X,Y]$, it follows that 
$$\{\{f,g\},s\}={\mathcal L}_{[X,Y]}s=[{\mathcal L}_X,{\mathcal L}_Y]s=\{f,\{g,s\}\}-\{g,\{f,s\}\}.$$
This clearly holds for any power $K^m$. 
\end{ex}

\begin{rmk} A holomorphic first-order operator $\bar D:{\mathcal O}(V)\rightarrow {\mathcal O}(V\otimes T)$ is globally defined as a vector bundle homomorphism $\alpha:J^1(V)\rightarrow V\otimes T$ where $J^1(V)$ is the bundle of holomorphic 1-jets of sections of $V$. It is an extension
$$0\rightarrow V\otimes T^*\rightarrow J^1(V)\rightarrow V\rightarrow 0$$
and its extension class in $H^1(M,\End V\otimes T^*)$, the Atiyah class, is the obstruction to splitting the sequence holomorphically. When $V$ is a line bundle, and $M$ is K\"ahler, this is the first Chern class in $H^{1,1}$.

The symbol $\sigma$ of $\bar D$ is the homomorphism $\alpha$ restricted to $V\otimes T^*\subset J^1(V)$, so
the existence of $\bar D$ means that $\sigma\in H^0(M,\Hom(V\otimes T^*,V\otimes T)$ lifts to a class $\alpha\in H^0(M,\Hom(J^1(V),V\otimes T)$. In the long exact cohomology sequence of the extension, this means that the map
$$\sigma: H^1(M,\End V\otimes T^*)\rightarrow H^1(M,\End V\otimes T)$$
applied to the Atiyah  class is zero. 

In the case of a line bundle this is the cup product we encountered in Section \ref{deform} applied to the first Chern class, so in particular we see that the existence of a Poisson module structure on the canonical bundle means that the image of $c_1(T)$ in $H^1(M,T)$ is zero. So, as in Section \ref{deform}, if $c_1(T)$ is represented by a K\"ahler form, Goto's theorem, to first order, keeps the complex structures $I_+,I_-$  in the same diffeomorphism class. 
\end{rmk}
 \vskip .25cm
 Just because the Lie derivative of a Hamiltonian vector field does not make the tangent bundle a Poisson module does not mean that it can't be one. If we take two  vector fields $X_1,X_2$ on $\CP^2$ then $\sigma=X_1\wedge X_2$ defines a Poisson structure. It is holomorphic symplectic where $\sigma$ is non-zero, which is away from a cubic curve $C$ -- the curve where $X_1\wedge X_2=0$ i.e. where $X_1$ and $X_2$ become linearly dependent.
 
 Here we can step back and view $\CP^2$ as a generalized complex manifold: away from $C$, $\sigma^{-1}$ defines a holomorphic section of the canonical bundle $K$ which we can write as a closed complex 2-form $B+i\omega$. The generalized complex structure here is a symplectic structure $\omega$ transformed by the B-field $B$. But on such a structure, a generalized holomorphic bundle is a flat vector bundle. Now $X_1$ and $X_2$ are linearly independent away from $C$ so we can try and define $\bar D$ on $T$ by making them covariant constant i.e. $\bar D X_1=\bar D X_2=0$. Then we need to show that this extends as a holomorphic differential operator -- the $\bar D^2=0$ condition is already satisfied on an open set and so holds everywhere. 
 
 Take  a local holomorphic basis $\partial/\partial z_1,\partial/\partial z_2$ for $T$ in a neighbourhood of a point of $C$, and then
 $$X_i=P_{ji}\frac{\partial}{\partial z_j}$$
 so
 \begin{equation}
 \sigma=X_1\wedge X_2=\det P\frac{\partial}{\partial z_1}\wedge \frac{\partial}{\partial z_2}.
 \label{x1x2}
 \end{equation}
A``connection matrix" for $\bar D$ relative to this basis  is given by a matrix $A$ of vector fields such that 
$$0=\bar DX_i=D(P_{ji}\frac{\partial}{\partial z_i})=\sigma(dP_{ji})\frac{\partial}{\partial z_j}+P_{ji}A_{kj}\frac{\partial}{\partial z_k}$$
which has solution
$$A=-\sigma(dP)P^{-1}=-\sigma(dP)\frac{\adj P}{\det P}.$$
From (\ref{x1x2}) this is 
$${\adj P}\left( \frac{\partial P}{\partial z_2}\frac{\partial}{\partial z_1}-\frac{\partial P}{\partial z_1}\frac{\partial}{\partial z_2}\right)$$
which is holomorphic and so $\bar D$ is well-defined.
\vskip .25cm
We can extend this argument to other rank 2 vector bundles $V$ with $\Lambda^2V\cong K^*$, so long as they have two sections $s_1,s_2$ to replace the vector fields $X_1,X_2$. A generic vector field on $\CP^2$ has three zeros: suppose $V$ has a section $s$ with $k$ simple zeros $x_1,\dots,x_k$ (and then the second Chern class $c_2(V)=k$). Then $s$ defines an exact  sequence of sheaves 
$$0\rightarrow {\mathcal O}\stackrel{s}\rightarrow {\mathcal O}(V)\rightarrow K^*\otimes {\mathcal I}
\rightarrow 0$$
where ${\mathcal I}$ is the ideal sheaf of the $k$ points. If $H^0(\CP^2,K^*\otimes {\mathcal I})\ne0$, in other words if there is a cubic curve $C$ passing through the $k$ points,  then from the exact cohomology sequence (and using $H^1(\CP^2,{\mathcal O})=0$) we can find a second section $s_2$ of $V$, and if $s_1=s$, $s_1\wedge s_2$ vanishes on the curve $C$, which defines a holomorphic Poisson structure. 

The {\it Serre construction} provides a means of constructing such bundles (see for example \cite{DK} Section 10.2.2). Away from the $k$ points we have an extension of line bundles
$$0\rightarrow {\mathcal O}\stackrel{s}\rightarrow {\mathcal O}(V)\rightarrow K^*
\rightarrow 0$$
which is described by a Dolbeault representative $\alpha\in \Omega^{0,1}(K)=\Omega^{2,1}$. 
It extends to an extension as above if it has a singularity at each of the points of the form 
$$\frac{1}{4r^4}dz_1\wedge dz_2\wedge (\bar z_2d\bar z_1-\bar z_1d\bar z_2).$$
In distributional terms $\bar\partial\alpha=\sum_i \lambda_i\delta_{x_i}=\beta$ a linear combination of delta functions of the points.

Such a sum  defines a class in $H^2(M,K)$. Since $H^2(M,K)$ is dual to $H^0(M,{\mathcal O})\cong \C$, this class is determined by evaluating it on the function $1$. But
$$\langle \beta, 1\rangle =\sum_i\lambda_i$$
so if the $\lambda_i$ sum to zero the class is zero and one can solve $\bar\partial \alpha=\beta$ for $\alpha$. 

Thus, given a cubic curve and a collection of $k$ points with non-zero scalars $\lambda_i$ whose sum is zero, we obtain a rank $2$ Poisson module with $c_2(V)=k$.

\begin{rmk} The Serre construction can also be used to generate co-Higgs bundles on $\CP^2$ with nilpotent Higgs field $\phi$. This time we require $\Lambda^2V\cong {\mathcal O}(1)$ and we want to solve $\bar\partial\alpha=\beta$ for a distribution defining a class in $H^2(\CP^2,{\mathcal O}(-1))$ which is dual to $H^0(\CP^2,{\mathcal O}(-2))=0$. Hence there is no constraint on the $\lambda_i$s. We obtain an extension
 $$0\rightarrow {\mathcal O}\stackrel{s}\rightarrow {\mathcal O}(V)\stackrel{\pi}\rightarrow {\mathcal O}(1)\otimes {\mathcal I}\rightarrow 0.$$
Choosing a section $w$ of $T(-1)$, $v\mapsto s(w\pi(v))$ defines $\phi \in H^0(\CP^2, \End V\otimes T)$ whose kernel and image lie in the trivial rank one subsheaf.   
\end{rmk}

\section{Holomorphic bundles and the B-field action}
\subsection{The B-field action}
On a complex manifold we can pull back a holomorphic vector bundle by a holomorphic diffeomorphism to get a new one, but in generalized geometry we have learned that the group $\Omega^2(M)_{cl}\rtimes \Diff(M)$ replaces the group of diffeomorphisms and in particular that a closed real $(1,1)$-form $B$ preserves the generalized complex structure determined by an ordinary complex structure.  We shall study next the effect of this action on generalized holomorphic bundles.

Recall that in this case a generalized holomorphic bundle is defined by an operator
  $$\bar Ds=\left(\frac{\partial s}{\partial \bar z_{i}}+A_{\bar i}s\right) d\bar z_{i}+\phi^{j} s \frac{\partial }{\partial z_{j}}$$
 where $(V,\phi)$ is a co-Higgs bundle.
 
 The transform of this by $B$ is then the operator
  $$\bar Ds=\left(\frac{\partial s}{\partial \bar z_{i}}+A_{\bar i}+\phi^{j}B_{j\bar i}s\right) d\bar z_{i}+\phi^{j} s \frac{\partial }{\partial z_{j}}$$

   More invariantly we write $i_{\phi}B\in \Omega^{0,1}(\End V)$ for the contraction of the Higgs field $\phi\in H^0(M,\End V\otimes T)$ with $B\in \Omega^{0,1}(T^*)$ and then we have a new holomorphic structure 
  $$\bar\partial_B=\bar\partial+i_{\phi}B$$
  on the $C^{\infty}$ bundle $V$. But the condition $\phi^2=0$ means that $i_{\phi}B$ commutes with $\phi$ and so $\phi$ is still holomorphic with respect to this new structure: $\bar\partial_B\phi=0.$
 
\begin{rmk} Note that if $U\subset V$ is a holomorphic subbundle with respect to $\bar\partial$ then  if it is also $\phi$-invariant, it is holomorphic with respect to $\bar\partial_B$. So stability is preserved by the B-field action.
\end{rmk}

 Now suppose  $B=\bar\partial \theta$ for  $\theta \in \Omega^{1,0}$ and define $\psi=i_{\phi}\theta$, a section of $\End V$. In coordinates  $\psi=\phi^{i}\theta_{i}$ which implies 
$(\bar\partial\psi)_{\bar i}=\phi^{j}B_{j\bar i}$. Then
$$[\psi,\bar\partial\psi]_{\bar i}=[\phi^{k}\theta_{k},\phi^{j}B_{j\bar i}]=[\phi^{k},\phi^{j}]\theta_{k}B_{j\bar i}=0$$
since $\phi^j$ and $\phi^k$ commute. This means  (unusually for a non-abelian gauge theory) that if we exponentiate to an automorphism of the bundle $V$ we have
 $$\exp (-\psi)\bar\partial\exp \psi=\bar\partial\psi +i_{\phi}B\psi.$$
 Thus if $B_1$ and $B_2$ represent the same Dolbeault cohomology class in $H^1(M,T^*)$, the two actions are related by an automorphism. Hence $H^1(M,T^*)$ acts on the moduli space of stable co-Higgs bundles. 
 
 \begin{ex} If we take the canonical Higgs bundle $V={\mathcal O}\oplus T$ and  $\phi(\lambda,X)=(X,0)$  as in Section \ref{coH}, then for $[B]\in H^1(M,T^*)$ the structure $\bar\partial_B$ defines a non-trivial extension 
 $$0\rightarrow {\mathcal O}\rightarrow V\rightarrow T\rightarrow 0$$
 which still has a canonical Higgs field. 
 \end{ex}

 We shall investigate this action in more detail next. 
 
\subsection{Spectral covers}
In the case of Higgs bundles, Simpson reinterpreted in \cite{S1} a {\it Higgs sheaf} on $M$ in terms of a sheaf on $\PP({\mathcal O}\oplus T^*)$ whose support is disjoint from the divisor at infinity.  This can be adapted immediately replacing $T^*$ by $T$.

In standard local coordinates $y_i,z_j$ on $TM$ given by the vector field  $y^i\partial/\partial z_i$, we pull back the rank $k$ bundle $V$ under the projection $p:TM\rightarrow M$ and define an action of $y^i$ by $\phi^i$. Since $[\phi^i,\phi^j]=0$ this defines a module structure over the commutative ring of functions polynomial in the fibre directions. 

More concretely, suppose in a  neighbourhood of a point some linear combination of the $\phi^i$, say $\phi^1$,  has distinct eigenvalues. Then since by the Cayley-Hamilton theorem $\phi^1$ satisfies its characteristic equation, on the support of the sheaf $y_1$ is an eigenvalue of $\phi^1$, and the kernel of $\phi^1-y_1$ defines a line bundle $U\subset p^*V$. Since all $\phi^i$ commute with $\phi^1$, $L$ is an eigenspace for $\phi^i$ with eigenvalue $y^i$. If the  $m$ characteristic equations of $\phi^i$ are generic, they define an $m$-dimensional  submanifold $S$ of $TM$ which is an $m$-fold covering of $M$ under $p$.  There will be points at which $\phi^1$ has coincident eigenvalues, but under suitable genericity conditions $S$ will still be smooth with a line bundle $L$. The action of $\phi$ on $U$ is
$$\phi\vert_U=\phi^i\vert_U\frac{\partial}{\partial z_i}=y_i\frac{\partial}{\partial z_i}$$
which is the tautological section of $p^*T$ on the total space of $T$.

The one-dimensional case of this, where $M=\CP^1$,  was much studied from the point of view of integrable systems before its important application to the moduli space of Higgs bundles, and the co-Higgs situation on $\CP^1$  is a particular case described, for example, in \cite{NJH5}. In this case, where $T={\mathcal O}(2)$,  $\phi$ is a holomorphic section of $\End V(2)$ and its characteristic equation is 
$$\det(\eta-\phi)=\eta^k+a_1\eta^{k-1}+\dots +a_k=0$$
where $a_i$ is a section of ${\mathcal O}(2i)$ on $\CP^1$. Interpreting $\eta=yd/dz$ as the tautological section of 
$p^*{\mathcal O}(2)$, this is the vanishing of a section of $p^*{\mathcal O}(2k)$ on the algebraic surface $T\CP^1$ and  it defines a spectral curve $S$ which, by the adjunction formula, has genus $g=(k-1)^2$  and is a branched covering of $\CP^1$ of degree $k$.

We reconstruct a co-Higgs bundle by taking the direct image $p_*L$ of a line bundle $L$ on $S$.  For any open set $U\subseteq \CP^1$, by definition
$$H^0(U,p_*L)=H^0(p^{-1}(U),L).$$
The sheaf $p_*L$ defines a rank $k$ vector bundle and the direct image of multiplication by the tautological section
$$\eta: H^0(p^{-1}(U),L)\rightarrow H^0(p^{-1}(U),L(2))$$
defines a Higgs field $\phi$. 

(Note that the line bundle $L$ is not quite the same as the eigenspace bundle $U$. The direct image gives a canonical evaluation map $p^*V\mapsto L$ so that $L^*\subset p^*V^*$ is the eigenspace bundle of the dual endomorphism $\phi^t$.)

The degree of the line bundle $L$ and that of the vector bundle $V$ are easily related -- the direct image  definition implies that $H^0(S,L\otimes p^*{\mathcal O}(n))\cong H^0(\CP^1, p_*L(n))$ and taking $n$ large, these are given by the Riemann-Roch formula. The result is 
$$\deg V=\deg L+k-k^2.$$
 \begin{exs}
 
 \noindent 1. Consider the example of $V={\mathcal O}\oplus {\mathcal O}(-1)$ in the previous lecture. Here $k=2$ and $\deg V=-1$, so $\deg L=1$. The curve $S$ has genus $(k-1)^2=1$ and is an elliptic curve ($y^2=c_0+c_1z+\dots +c_4z^4$). The line bundle $L$ has degree one and hence has a unique section which vanishes at a single point, which is $\eta(0)$ in our description of the moduli space.
 
 \noindent 2. If $\deg L= g-1= k^2-2k$ then $\deg V=-k$, so $V(1)$ has degree zero. Now a vector bundle $E$ on $\CP^1$ is trivial if it has degree zero and $H^0(\CP^1,E(-1))=0$, so $V(1)$ is trivial if $0=H^0(\CP^1,V)=H^0(S,L)$, which is if the divisor class of $L$ does not lie on the theta divisor of $S$. In this case a co-Higgs bundle consists of a $k\times k$ matrix whose entries are sections of ${\mathcal O}(2)$.
 \vskip .25cm
 Now consider the B-field action from the point of view of the spectral cover. Since $\phi$ is only changed by conjugation, the spectral cover is unchanged -- it is only the holomorphic structure on the line bundle which can change. But the change in the holomorphic structure on $V$ was 
 $$\bar\partial \mapsto \bar\partial +i_\phi B$$
 and on $U\subset V$ $\phi$ acts via the tautological section $\eta$ of $p^*TM$, so we are changing the holomorphic structure of $U$ by
 $$\bar\partial \mapsto \bar\partial +i_\eta B.$$
 In other words, we have $[B]\in H^1(M,T^*)$ which we pull back to $p^*[B]\in H^1(TM,p^*T)$ then contract with $\eta\in H^0(TM,p^*T)$ to get the class
 $$\eta p^*[B]\in H^1(TM,{\mathcal O}).$$
 Exponentiating to $H^1(TM,{\mathcal O}^*)$ defines a line bundle $L_B$. Restricting to $S$ the B-field action is  $U\mapsto U\otimes L_B$.
 \end{exs}
 
 Let us look at this action in the two examples above. Since $H^{1,1}(\CP^1)$ is one-dimensional a real closed $(1,1)$ form is cohomologous to a  multiple of 
 $$\frac{idz\wedge d\bar z}{(1+z \bar z)^2}$$ 
This form  integrates to $2\pi$ over $\CP^1$. 

Pulling back to $T\CP^1$ and contracting with $yd/dz$ we obtain the class in $H^1(S,{\mathcal O})$ represented by
 $$\frac{iy d\bar z}{(1+z \bar z)^2}.$$ 
 
 \noindent 1. In the first example,  $S$ is an elliptic curve and a point of the moduli space is defined by a point $x$ on this curve, so  tensoring with a line bundle $L_B$ is a translation. The  non-vanishing 1-form $dz/y$ is equal to $du$ in the uniformization and then two points $x,x'$ are related by a translation $u\mapsto u+a$ if
 $$\int_x^{x'}\frac{dz}{y}=a$$
 modulo periods. 
 
 On the other hand our class in $H^1(S,{\mathcal O})$  pairs with $dz/y\in H^0(S, K)$ by integration:
$$\int_S \frac{idz\wedge d\bar z}{(1+z \bar z)^2}=4\pi$$
so this determines the translation.

\noindent 2. In the second example, since the bundle is trivial we may write the Higgs field in $\End \C^k\otimes H^0(\CP^1,{\mathcal O}(2))$  as a matrix with entries quadratic polynomials in $z$. Write it thus:
$$\phi=(T_1+iT_2)+2iT_3 z+(T_1-iT_2)z^2.$$
 Then, as derived in \cite{NJH5}, tensoring by $L_B$ is integrating to time $t=1$ the system of nonlinear differential equations called {\it Nahm's equations}.
 $$\frac{dT_1}{dt}=[T_2,T_3],\quad \frac{dT_2}{dt}=[T_3,T_1],\quad \frac{dT_3}{dt}=[T_1,T_2].$$
These equations arise in the study of non-abelian monopoles and are dimensional reductions of the self-dual Yang-Mills equations. 

From these examples it is clear that the B-field action can be highly non-trivial. What it also shows is that the action on the moduli space can be quite badly behaved, for the Nahm flow could be an irrational flow on the Jacobian of the spectral curve. 

\subsection{Twisted bundles and gerbes}
Now suppose we replace the generalized complex structure on $T\oplus T^*$ by a twisted version on the bundle $E$ defined by a 1-cocycle $B_{\alpha\beta}$ of closed real $(1,1)$-forms. What is a generalized holomorphic bundle now? The general definition is the same -- a vector bundle $V$ with a differential operator $\bar D$ but we want to understand it in more concrete terms. 

If we think of $E$ as obtained by patching together copies of $T\oplus T^*$ then over each open set $U_{\alpha}$, $V$ has the structure of a co-Higgs bundle -- a holomorphic structure $A_{\alpha}$ and a Higgs field $\phi_\alpha$. On the intersection $U_{\alpha}\cap U_{\beta}$ these are related by the B-field action of $B_{\alpha\beta}$:

 \begin{equation}
 (A_{\beta})_{\bar i}=(A_{\alpha})_{\bar i}+\phi^{j}(B_{\alpha\beta})_{j\bar i}, \qquad (\phi_{\beta})^{j} = (\phi_{\alpha})^{j}.
 \label{af}
 \end{equation}

Consider first the case of $V=L$ a line bundle. Then, because $\End V$ is holomorphically trivial for all  of the local holomorphic structures, $\phi$ is a global holomorphic vector field $X$. So consider the $(0,1)$ form 
$$A_{\alpha\beta}=i_XB_{\alpha\beta}.$$
The $(1,1)$ form $B_{\alpha\beta}$ is closed so $\bar\partial B_{\alpha\beta}=0$ and $X$ is holomorphic so that $\bar\partial A_{\alpha\beta}=0$. Locally write $A_{\alpha\beta}=\bar\partial f_{\alpha\beta}$, then, since $B_{\alpha\beta}$ is a cocycle, on threefold intersections $f_{\alpha\beta}+f_{\beta\gamma}+f_{\gamma\alpha}$ is holomorphic. Write 
$$g_{\alpha\beta\gamma}=\exp 2\pi i (f_{\alpha\beta}+f_{\beta\gamma}+f_{\gamma\alpha})$$
then this defines a holomorphic gerbe.

But the local holomorphic structure on $L$ is defined by a $\bar\partial$-closed form $A_{\alpha}$, so writing $A_{\alpha}=\bar\partial h_{\alpha}$ we have from (\ref{af}) that $k_{\alpha\beta}=f_{\alpha\beta}+h_{\alpha}-h_{\beta}$ is holomorphic and moreover 
$$g_{\alpha\beta\gamma}=\exp 2\pi i (k_{\alpha\beta}+k_{\beta\gamma}+k_{\gamma\alpha}).$$
This is a {\it holomorphic trivialization} of the gerbe, or as is sometimes said, a line bundle over the gerbe. The ratio of any two trivializations (i.e. writing $g_{\alpha\beta\gamma}$ as a coboundary) is a cocycle which defines the transition functions for a holomorphic line bundle. In the untwisted case a generalized holomorphic bundle was  just  a line bundle and a vector field;  here any two {\it differ} by such an object. 

In more invariant terms we have taken the class in $H^2(M,T^*)$ defined by the $(1,2)$ component of the 3-form $H$, and contracted with the vector field $X\in H^0(M,T)$ to get a class in $H^2(M,{\mathcal O})$. Exponentiating gives us an element in $H^2(M,{\mathcal O}^*)$ which is the equivalence class of the holomorphic gerbe defined by $g_{\alpha\beta\gamma}$.  The existence of a trivialization of the gerbe is the statement that this class is zero.

\vskip .25cm
Now consider the general case: over each $U_{\alpha}$ we can consider the spectral cover in $TM$. This is defined by characteristic polynomials of components of $\phi$. The $C^{\infty}$ transition functions for the vector bundle $V$ conjugate $\phi$ and so leave these polynomials invariant. It follows that the local spectral covers fit together into a global spectral cover $S\subset TM$. The eigenspace bundle $U$ however, only has local holomorphic structures. But $\phi$ acts on $U$ via the tautological section $\eta$ of $p^*TM$, and so we are in a parallel situation to  the one we just considered: a gerbe on $TM$ defined by the cocycle
$$A_{\alpha\beta}=i_\eta p^*B_{\alpha\beta}.$$
In the untwisted case, a co-Higgs bundle was determined by a line bundle on the spectral cover, in this case it is a trivialization of the gerbe.
\vskip .25cm
The language of gerbes is convenient to describe things on the spectral cover, but a $C^{\infty}$ bundle $V$ with local holomorphic structures is not readily adaptable to  conventional algebraic geometric language on $M$ itself.  As far as generalized geometry is concerned we have $\bar D$, but it is still useful to rephrase the structure in more conventional language.  For that purpose, we can split the extension $E$ and work with $T\oplus T^*$ and the Courant bracket twisted with a 3-form $H$. 

The generalized Dolbeault complex is now $\bar D=\bar\partial_A-H^{1,2}+\phi$ where
$$\bar\partial_A: \Omega^{0,q}(V\otimes \Lambda^pT)\rightarrow \Omega^{0,q+1}(V\otimes \Lambda^pT),$$
 the 3--form $H^{1,2}$  acts by contraction in the $(1,0)$ entry,
 $$H^{1,2}:\Omega^{0,q}(V\otimes \Lambda^pT)\rightarrow \Omega^{0,q+2}(V\otimes \Lambda^{p-1}T)$$
 and the Higgs field acts like this
$$\phi:  \Omega^{0,q}(V\otimes \Lambda^pT)\rightarrow \Omega^{0,q}(V\otimes \Lambda^{p+1}T).$$
The condition $\bar D^2=0$ now becomes 
$$\bar\partial_A^2=i_{\phi}H,\quad \bar\partial_A\phi=0,\quad \phi^2=0.$$

This shape of structure has appeared in the literature. For example, replacing $V$ by $\End V$ (and thereby getting a complex which governs the deformation theory of a generalized holomorphic bundle), we obtain a {\it curved differential graded algebra} -- an  algebra with derivation where $d^2a=[c,a]$ and $dc=0$. This is an identifiable concept, but nevertheless, packaged in the language of generalized geometry it becomes quite natural.

\vskip 1cm
 Mathematical Institute, 24-29 St Giles, Oxford OX1 3LB, UK
 
 hitchin@maths.ox.ac.uk

 \end{document}